\documentclass[12pt]{article}

\usepackage[margin=1in]{geometry}

\usepackage{amsmath,amssymb,amsfonts}
\usepackage{mathtools}
\usepackage{bm}
\usepackage{graphicx}
\usepackage{booktabs}
\usepackage{siunitx}
\usepackage{enumitem}
\usepackage{multirow}
\usepackage{float}
\usepackage{hyperref}
\usepackage{xcolor}
\usepackage{caption}
\usepackage{subcaption}

\usepackage{algorithm}
\usepackage{algorithmicx}
\usepackage{algpseudocode}

\usepackage{tikz}
\usetikzlibrary{calc,positioning,shapes.geometric,arrows.meta}
\usepackage{tikz-3dplot}

\usepackage{listings}

\numberwithin{equation}{section}

\newcommand{\dt}{\Delta t}
\newcommand{\PsiLS}{\Phi}
\newcommand{\TT}{\mathcal{T}}
\newcommand{\Kcut}{\mathcal{K}_{\mathrm{cut}}}
\newcommand{\levelN}{\mathrm{levelN}}
\newcommand{\TTn}{\{\TT^0_n,\TT^1_n,\dots,\TT^{\levelN}_n\}}
\newcommand{\TTnp}{\{\TT^0_{n+1},\TT^1_{n+1},\dots,\TT^{\levelN}_{n+1}\}}

\newcommand{\eugenio}[1]{#1}

\definecolor{ForestGreen}{RGB}{34,139,34}

\title{An Adaptive Finite Element Method for Marker-Driven Level-Set Transport on Hierarchical Meshes}

\author{
Eugenio Aulisa\thanks{Department of Mathematics and Statistics, Texas Tech University, Lubbock, TX, USA}
\and
Samuele Baldini\thanks{Department of Industrial Engineering, University of Bologna, Bologna, Italy}
\and
Giacomo Barbi\footnotemark[1]
\and
Andrea Chierici\footnotemark[1]
\and
Sandro Manservisi\footnotemark[2]
}

\date{}

\begin{document}

\maketitle

\noindent\textbf{Accepted manuscript.} The final authenticated version is available at:\\
\url{https://doi.org/10.1016/j.jcp.2026.115229}

\begin{abstract}
\eugenio{This work presents a new adaptive finite-element framework for level-set transport, achieving high accuracy in kinematic interface transport problems relevant to interface-capturing methods for two-phase flows while reducing computational cost.}
This framework accommodates dynamic refinement and coarsening, and is compatible with standard finite-element data structures. The algorithm is applicable to both structured and unstructured discretizations in two and three dimensions. The method uses a tree-based hierarchical mesh with dynamic local refinement that tracks the evolving interface, concentrating resolution in a narrow band around the zero level set while preserving a coarse discretization elsewhere. The level-set field is updated through marker transport on adaptively refined meshes, enabling interface prediction and the construction of an evolving adaptive hierarchy. Concurrently, backward characteristic tracing provides accurate evaluation of the advected level-set field. Both operations employ an efficient multilevel marker and point-location algorithm to identify containing elements across refinement levels. Numerical experiments on two- and three-dimensional structured and unstructured grids, including quadrilateral, simplicial, wedge, and hexahedral meshes\eugenio{, subjected to severe interface deformations,} demonstrate that the adaptive strategy achieves accuracy comparable to uniform discretizations at a substantially lower cost, \eugenio{while maintaining good conservation properties under aggressive adaptivity. Consequently, the proposed approach provides an efficient and flexible framework that can be naturally integrated into level-set-based multiphase flow solvers.}
 
\end{abstract}


\section{Introduction}
\subsection{Background}

The simulation of two-phase and, more generally, multiphase flows requires an accurate and robust description 
of moving interfaces, presenting significant numerical challenges when handling large deformations, 
topological changes, and coupled physical phenomena. 
\eugenio{In the present work, we focus on the kinematic transport component of such interface-capturing formulations.}
The numerical challenge is to advance these interfaces in the presence of complex velocity fields, 
inducing strong deformations \eugenio{that can lead to potential topological
changes}. To this end,
two different techniques for transporting the interface have been developed: explicit interface-tracking 
methods, 
which maintain a sharp, discrete representation of the interface geometry, and implicit interface-capturing methods, which embed the interface within a fixed grid using auxiliary scalar fields. The 
choice between these approaches involves inherent trade-offs between geometric precision, conservation 
properties, and robustness to topological changes.

Among these approaches, the level-set method \cite{OsherSethian1988,sussman1994level} addresses this 
issue by transporting a smooth (Lipschitz-continuous) scalar field $\PsiLS(x,t)$, whose zero isoline or 
isosurface represents the moving interface. This technique avoids any explicit treatment of the 
interface, thereby facilitating complex interface dynamics and topological changes.
\eugenio{In particular, topological changes are handled implicitly through the Eulerian evolution of the level-set field rather than through explicit interface reconstruction.}
Moreover, geometric 
quantities of interest, such as normals and curvature, can be extracted directly from the level-set 
field, making this formulation compatible with a wide range of discretization and differentiation 
schemes. This formulation is compatible with both finite-volume/finite-difference and finite-element 
discretizations, as geometric quantities are derived directly from the level-set field.
Despite its elegance and versatility, the level-set method suffers from several well-documented 
drawbacks, most notably mass loss and degradation of the level-set field, both of which can introduce 
significant errors during advection. Mass loss typically arises from discretization inaccuracies and the 
numerical diffusion required to stabilize the advection process, undermining volume conservation. 
Meanwhile, degradation refers to the departure of the level-set function from a signed-distance profile due to numerical advection, necessitating periodic reinitialization to restore accuracy 
\cite{osher2004level,min2010reinitializing}.

Historically, these difficulties have been addressed using several strategies, including high-order 
Hamilton–Jacobi solvers (ENO/WENO schemes with TVD Runge–Kutta integration \cite{JiangPeng2000}), narrow-
band evolution \cite{AdalsteinssonSethian1995}, reinitialization of $\PsiLS$ toward a signed-distance 
function \cite{sussman1994level}, and hybridizations with particles or VOF techniques to improve mass 
conservation \cite{Enright2002,sussman2000coupled}.
\eugenio{In many particle-assisted formulations, particles contribute directly to interface correction or reconstruction. In contrast, the markers introduced in the present work are not used to reconstruct the interface itself.}
For a broader overview of level-set formulations, 
transport schemes, reinitialization strategies, and their applications to multiphase flow, we refer to 
the review by Gibou et al. \cite{gibou2018review}.
Since the level-set method extends the interface motion to a volumetric field to simplify its advection, 
the region where high accuracy is most critical is the neighborhood of the zero level-set. Moreover, two-
phase phenomena may occur at very small scales, demanding fine spatial detail and thus incurring high 
computational cost. For this reason, adaptive mesh refinement (AMR) aligns naturally with level-set 
methods, since it concentrates resolution at the interface. The objective is to use only 
the computational effort required to capture the relevant physical details. Adaptive strategies are 
typically divided into h-adaptivity and p-adaptivity. In h-adaptivity, the mesh is locally refined by 
increasing the number of elements, which reduces the local mesh size $h$, a principle underlying most AMR 
frameworks \cite{berger1984amr1,berger1989amr2}. In p-adaptivity, the mesh remains fixed while the local 
polynomial degree of the basis functions is increased to improve accuracy, a strategy grounded in the 
classical p-version finite element theory \cite{Babuska1981}. Recent work has applied p- and hp-adaptive 
high-order discretizations to multiphase flows \cite{ntoukas2022entropy,mossier2023efficient}. 

In this work, we focus exclusively on h-adaptivity. Within its context, one major group of methods is 
block-structured Cartesian AMR, where the mesh is organized as a hierarchy of nested Cartesian grids 
\cite{sussman1999amr}. A recent formulation \cite{zeng2023amr} exemplifies current developments in this 
class, incorporating consistent mass–momentum transport and enabling stable simulations at high density 
ratios and Reynolds numbers. A second group is tree-based AMR, in which refinement proceeds by 
recursively subdividing individual cells to form a hierarchical tree structure 
\cite{khokhlov1998amr,Popinet2003}. 
The use of quadtree- and octree-based grids for level-set computations was first introduced in 
\cite{strain1999amr} and subsequently extended in several works. A second-order level-set method on such 
grids was presented in \cite{MinGibou2007}. The quadtree/octree framework of Popinet 
\cite{popinet2015quadtree} further advanced this class of methods by providing a robust adaptive solver 
for multiphase incompressible flow. While these approaches remain tied to structured Cartesian trees, in 
this work we combine semi-Lagrangian advection with a tree hierarchy on unstructured finite element 
meshes.

More recent advances include parallel implementations of Cartesian tree-based AMR 
\cite{burstedde2011p4est,mirzade2016amr}, improved interface-capturing strategies on adaptive quadtrees 
\cite{hergibo2024quadtree}, and updated octree-based formulations for fluid dynamics 
\cite{yang2025octree}. These methods rely on a Cartesian-cell refinement paradigm driven by uniform 
quadtree or octree subdivision. In contrast, our approach embeds the tree hierarchy directly within an 
unstructured finite element mesh and combines it with semi-Lagrangian transport. This allows adaptive 
refinement without being constrained to a Cartesian background grid, enabling resolution to be increased 
freely based on geometric or solution features rather than a fixed splitting pattern.

This distinction also extends to classical adaptive finite element methods (AFEM), where elements of 
general shape are refined under mesh conformity constraints rather than by uniform Cartesian subdivision 
\cite{ngo2020multi}. While Cartesian AMR benefits from regularity and simple indexing, numerical 
simulation of multiphase flows on unstructured meshes introduces additional challenges due to complex 
connectivity, irregular geometries, and the absence of a uniform refinement template. However, adaptive 
refinement on tetrahedral and other unstructured grids offers the geometric flexibility required for 
complex, body-fitted discretizations.
Early algorithms for local refinement of tetrahedral meshes were introduced in \cite{zhang2009parallel}, 
and more recent unstructured tetrahedral AMR strategies appear in \cite{morgan2017tetra}. Related octree-
based finite element frameworks and multigrid solvers on hierarchical meshes have been proposed 
\cite{sampath2010parallel}, yet level-set transport and reinitialization on such adaptive FEM trees 
remain largely unaddressed. Recent adaptive FEM frameworks have begun to integrate level-set formulations 
for multiphase flows with surface tension on dynamically refined triangular and tetrahedral meshes 
\cite{ngo2020multi}, permitting localized resolution near evolving interfaces while preserving the 
geometric adaptability of unstructured grids.
\eugenio{Our work contributes to this area} by introducing a tree-based adaptive FEM framework that supports efficient, 
mesh-aware semi-Lagrangian level-set transport and reinitialization, extending the flexibility of 
unstructured adaptation to interface-capturing multiphase flow simulations.

For a prescribed velocity field ${v}$, the evolution equation for $\PsiLS$ can be written as a linear 
advection (Hamilton--Jacobi) problem of the form
\[
\partial_t \PsiLS + {v}\cdot\nabla \PsiLS = 0.
\]
Accuracy and stability are often enhanced by high-order upwind spatial
discretizations, such as ENO/WENO schemes, coupled with strong-stability-preserving (SSP) or total-
variation-diminishing (TVD) time integrators~\cite{JiangPeng2000,Shu2009}. 
An alternative approach is provided by semi-Lagrangian (SL) schemes, which update the solution by 
integrating characteristics backward in time and interpolating at the departure points. This relaxes the 
usual CFL constraint and enables comparatively large time steps while retaining good stability 
properties~\cite{StaniforthCote1991}. Although SL schemes do not enforce mass conservation, they avoid 
the additional mass errors introduced by strongly diffusive upwind stabilizations. Its main drawbacks are 
the need for accurate and efficient interpolation and the associated computational overhead.
Several works have explored SL-based approaches for interface transport, including semi-Lagrangian 
contouring for level-set advection \cite{bargteil2006semi} and improved particle level-set formulations 
designed to reduce mass-loss errors \cite{kurioka2025improved}.
\eugenio{These approaches typically employ severe deformation benchmarks to assess robustness under thin filament formation and near-topological changes \cite{Rider1998, Enright2002, Wang2009, Cervone2009, Ramanuj2019}.}

To prevent the progressive distortion of the level-set field induced by advection, $\PsiLS$ is 
periodically reinitialized to recover a signed-distance profile while preserving the zero level-set.
Classical techniques fall into three broad families: direct distance computation, fast-marching or fast-
sweeping methods, and PDE-based formulations.
The direct approach computes the signed-distance function by evaluating, at each node, the minimum 
Euclidean distance to the interface \cite{shakoor2025reinit}. Fast-marching \cite{sethian2001reinit} and 
fast-sweeping \cite{zhao2005reinit} methods solve the Eikonal equation $|\nabla d|=1$ with $d=0$ on the 
interface to obtain a signed-distance field $\PsiLS=\operatorname{sign}(\PsiLS_0),d$. These methods are 
highly efficient but typically first-order near the interface. The PDE-based formulation 
\cite{sussman1994level} computes a signed-distance map as the steady state of
\[
\frac{\partial \PsiLS}{\partial \tau} + \operatorname{sign}(\PsiLS_0)\,\bigl(|\nabla \PsiLS| - 1\bigr) = 
0,
\]
where $\tau$ is a pseudo-time variable and $\PsiLS_0$ is the initial level-set function. This equation is 
advanced in pseudo-time until a steady state is reached, at which point $\PsiLS$ approximates a signed-
distance function. This approach can achieve high accuracy when discretized with high-order Hamilton–
Jacobi schemes. Variants that improve curvature accuracy and mass preservation include conservative and 
smooth reinitialization strategies \cite{Olsson2007}.
A fourth line of work adopts a geometric viewpoint based on closest-point projection. Here the signed-
distance value at a node is reconstructed by projecting onto the interface and measuring its Euclidean 
distance. Early formulations include \cite{chopp2001reinit}, which combines projection near the interface 
with fast marching elsewhere, and \cite{anumolu2013reinit}, which follows a similar hybrid strategy in 
finite differences. Extensions improving robustness through kink detection appear in 
\cite{henri2022reinit}. Closest-point ideas have also been introduced in the finite-element setting, such 
as the exact projection method for linear elements in \cite{parolini2007reinit}. These geometric 
approaches avoid solving a nonlinear Hamilton–Jacobi problem but require accurate nearest-point searches 
and careful treatment of the interface zone.

\subsection{A new mixed
marker and level-set front-tracking algorithm}

In this work, \eugenio{we present a level-set solver on marker-driven adaptively refined hierarchical meshes}, combining a tree-based 
adaptive mesh refinement (AMR) strategy with a semi-Lagrangian transport scheme.
\eugenio{The proposed method is intended as a highly adaptive transport framework that can be integrated into fully coupled multiphase flow solvers.}
This approach extends 
the mixed marker and level-set front-tracking framework introduced in \cite{aulisa2025mixed} and is 
\eugenio{designed to improve the preservation of fine-scale interface structures, eliminate refinement constraints typical of Eulerian formulations}, and enable efficient simulation of complex geometries.

The level-set function is discretized within a finite-element framework on hierarchical quadtree (2D) or 
octree (3D) meshes. Unlike single-mesh AMR formulations, we maintain an explicit multilevel hierarchy
$
 \TT^{0}\subset \TT^{1}\subset\cdots\subset \TT^{\levelN},
$ with stored parent–child relationships, and treat the entire hierarchy as a unified data structure 
throughout each time step. Refinement is dynamically concentrated within a narrow band surrounding the 
zero level-set, with the level-set field stored exclusively on the finest mesh level.

A central feature of our method is the construction of refinement in reference (canonical) coordinates. 
For quadrilateral and hexahedral elements, this corresponds to uniform Cartesian subdivision in reference 
space, while for triangular, tetrahedral, and wedge elements, refinement is defined via fixed templates 
on the master element. In all cases, the mapping to physical space ensures the underlying mesh remains 
unstructured. This preserves the simplicity and efficiency of Cartesian tree-based data structures, while 
greatly simplifying advection, interpolation, reinitialization, and point-location operations. This 
reference-space refinement strategy extends naturally to simplex elements, since subdivision occurs in 
the master element rather than through geometric splitting in physical coordinates.

Our tree-based refinement fundamentally differs from classical finite-element AMR based on element 
bisection or longest-edge subdivision. Conventional FEM AMR enforces mesh conformity constraints and 
generates refinement patterns dictated by mesh topology, complicating the dynamic redistribution of 
degrees of freedom and increasing the cost of characteristic-based transport. In contrast, our reference-
tree AMR produces a uniform hierarchical structure independent of element shape. Because refinement is 
performed in reference coordinates, hanging-node configurations follow a predictable hierarchy, and no 
mesh-conformity propagation is required. This significantly simplifies dynamic adaptivity and reduces 
refinement costs compared to traditional FEM strategies.

A second key contribution is the use of semi-Lagrangian advection to decouple the transport step from the 
refinement pattern and update the level-set field directly on the adapted mesh. Each time step proceeds 
in two phases: forward marker advection with mesh reconstruction and update.
First marker points extracted from cut elements are advected forward in time and used to reconstruct a new 
multilevel hierarchy.
\eugenio{The markers are used exclusively for geometric prediction and adaptive refinement and do not define interface topology or connectivity.}
Then, the level-set field on the new finest level is updated by tracing 
characteristics backward into the previous hierarchy and interpolating $\PsiLS$ at the departure points.

This forward–rebuild/backward–update structure eliminates the need for compatibility between successive 
adaptive meshes and avoids projecting the level-set field through intermediate coarse representations.
The hierarchical tree organization further accelerates point-location queries along characteristics via a 
coarse-to-fine push-down mechanism based on stored parent–child maps, making semi-Lagrangian advection 
particularly efficient in this multilevel AMR context.
Finally, we introduce a closest-point reinitialization algorithm tailored to the tree hierarchy, which 
reconstructs a signed-distance field without solving a nonlinear Hamilton–Jacobi equation. Its 
performance benefits directly from the structured refinement pattern near the interface.

While the mixed marker and level-set framework of \cite{aulisa2025mixed} already enabled two-dimensional 
dynamic simulations and benchmark flows, the present study focuses exclusively on kinematic interface 
transport and introduces a multilevel input/output formulation well-suited to highly adaptive semi-
Lagrangian updates.
\eugenio{Accordingly, the numerical validation presented in this work is restricted to prescribed kinematic velocity fields and severe interface deformation benchmarks.}
The multilevel mesh hierarchy developed here is designed for reuse in a fully 
dynamic formulation based on the incompressible Navier–Stokes equations. In that setting, the same 
hierarchy provides the natural data structure for multilevel and multigrid solvers on h-refined finite-
element spaces with arbitrary hanging-node configurations, as constructed in \cite{aulisa2019construction}. 
The extension to fully dynamic two-phase flow will therefore integrate the present adaptive semi-
Lagrangian transport framework with Navier–Stokes solvers based on multilevel algorithms, combined with 
the stabilization techniques introduced in \cite{gammanpila2025stabilized} for consistent treatment of 
interfacial stresses and discontinuous material properties.

\section{Notations on multilevel mesh construction} \label{sec:notation}

Let $\Omega $ be the computational domain and $ \TT^{0}(\Omega) $ be its coarse triangulation. We consider the adaptive multilevel finite-element mesh hierarchy
\[
\TT^{0}\subset \TT^{1}\subset \cdots \subset \TT^{\levelN},
\]
built from a uniform level-$0$ mesh by iterative quadtree (2D) / octree (3D) refinement subject to a one-level (2:1) grading constraint.
All levels $\TT^\ell$ are stored together with explicit father-children relations and are treated as a single multilevel data structure throughout.

\paragraph{Element types and reference domains}
The following higher-order finite elements are employed:
\begin{itemize}[leftmargin=*]
  \item Quad9: $\hat K = [-1,1]^2$,
  \item Tri7: $\hat K = \mathrm{conv}\{(0,0),(1,0),(0,1)\}$,
  \item Hex27: $\hat K = [-1,1]^3$,
  \item Tet15: $\hat K = \mathrm{conv}\{(0,0,0),(1,0,0),(0,1,0),(0,0,1)\}$,
  \item Wedge21: $\hat K = \hat K_{\mathrm{tri}}\times[-1,1]$.
\end{itemize}
For each element family, a canonical reference element, $\hat K$, and a consistent node ordering are taken as described in Figures~\ref{fig:quad9children}-\ref{fig:tet21children}.

\paragraph{Remark on the choice of higher-order element families and FEM pairing}
The adoption of Tri7, Tet15, and Wedge21 elements is dictated by consistency with the finite-element technology implemented in \textsc{FEMuS}~\cite{femus-github}, rather than by requirements intrinsic to the semi-Lagrangian formulation.
Standard alternatives such as Tri6, Tet10, and Wedge18 could in principle be used.

These enriched elements were introduced in \textsc{FEMuS} to obtain stable discretizations of the incompressible Navier-Stokes equations with biquadratic or triquadratic velocity and discontinuous linear pressure.
While the $Q_2$-$P_1^{\mathrm{DG}}$ pairing is stable on quadrilateral and hexahedral meshes, analogous simplicial and prismatic pairings based on Tri6, Tet10, or Wedge18 were observed to be unstable in practice within the \textsc{FEMuS} framework.

Stability is recovered by enriching the velocity space with interior degrees of freedom: a barycenter node in 2D (Tri7) and face-barycenter and interior nodes in 3D (Tet15 and Wedge21). Although no general inf-sup result is claimed, extensive numerical evidence demonstrates robust behavior on simplicial, prismatic, and hybrid meshes~\cite{aulisa2018monolithic,calandrini2019fluid,calandrini2020field}.

Accordingly, these element families are retained here for compatibility with Navier-Stokes solvers in \textsc{FEMuS}.
Their barycenter-enriched structure also supports consistent refinement templates and robust mappings under adaptive multilevel refinement, which are essential for reliable characteristic tracing and interpolation.



\begin{figure}
\centering

\begin{minipage}[t]{0.48\textwidth}
\centering
\begin{tikzpicture}[
  scale=1.5, 
  pt/.style={circle, fill=black, inner sep=1.1pt},
  lbl/.style={font=\scriptsize, inner sep=1pt},
  ed/.style={line width=0.6pt},
  dotline/.style={line width=0.7pt, draw=gray!65, dotted}
]
\coordinate (n0) at (-1,-1);
\coordinate (n1) at ( 1,-1);
\coordinate (n2) at ( 1, 1);
\coordinate (n3) at (-1, 1);
\coordinate (n4) at ( 0,-1);
\coordinate (n5) at ( 1, 0);
\coordinate (n6) at ( 0, 1);
\coordinate (n7) at (-1, 0);
\coordinate (n8) at (0,0);

\draw[ed] (n0)--(n1)--(n2)--(n3)--cycle;
\draw[dotline] (n7) -- (n5);
\draw[dotline] (n4) -- (n6);

\foreach \i in {0,1,2,3,4,5,6,7,8}{
  \node[pt] at (n\i) {};
  \node[lbl, anchor=west] at ($(n\i)+(0.03,0.02)$) {\i};
}

\node[anchor=west] at (1.35,0.2) {%
\scriptsize
\setlength{\tabcolsep}{3pt}
\renewcommand{\arraystretch}{1.05}
\begin{tabular}{c|cccc}
\toprule
child & \multicolumn{4}{c}{vertices} \\
\midrule
0 & 0 & 4 & 8 & 7 \\
1 & 4 & 1 & 5 & 8 \\
2 & 8 & 5 & 2 & 6 \\
3 & 7 & 8 & 6 & 3 \\
\bottomrule
\end{tabular}
};
\end{tikzpicture}

\captionof{figure}{Quad9 nodes and refinement.}
\label{fig:quad9children}
\end{minipage}
\hfill
\begin{minipage}[t]{0.48\textwidth}
\centering
\begin{tikzpicture}[
  scale=3.0,
  pt/.style={circle, fill=black, inner sep=1.1pt},
  ptgray/.style={circle, fill=gray!65, inner sep=1.1pt},
  lbl/.style={font=\scriptsize, inner sep=1pt},
  lblgray/.style={font=\scriptsize\color{gray!65}, inner sep=1pt},
  ed/.style={line width=0.6pt},
  dotline/.style={line width=0.7pt, draw=gray!65, dotted}
]
\coordinate (n0) at (0,0);
\coordinate (n1) at (1,0);
\coordinate (n2) at (0,1);
\coordinate (n3) at (0.5,0);
\coordinate (n4) at (0.5,0.5);
\coordinate (n5) at (0,0.5);
\coordinate (n6) at (0.333333,0.333333);

\draw[ed] (n0)--(n1)--(n2)--cycle;

\draw[dotline] (n3)--(n5);
\draw[dotline] (n3)--(n4);
\draw[dotline] (n5)--(n4);

\foreach \i in {0,1,2,3,4,5}{
  \node[pt] at (n\i) {};
  \node[lbl, anchor=west] at ($(n\i)+(0.02,0.015)$) {\i};
}
\node[ptgray] at (n6) {};
\node[lblgray, anchor=west] at ($(n6)+(0.02,0.015)$) {6};

\node[anchor=west] at (1.05,0.65) {%
\scriptsize
\setlength{\tabcolsep}{3pt}
\renewcommand{\arraystretch}{1.05}
\begin{tabular}{c|ccc}
\toprule
child & \multicolumn{3}{c}{vertices} \\
\midrule
0 & 0 & 3 & 5 \\
1 & 3 & 1 & 4 \\
2 & 5 & 4 & 2 \\
3 & 4 & 5 & 3 \\
\bottomrule
\end{tabular}
};
\end{tikzpicture}

\captionof{figure}{Tri7 nodes and refinement.}
\label{fig:tri7children}
\end{minipage}

\vspace{8mm}

\begin{minipage}[t]{0.92\textwidth}
\centering
\begin{tikzpicture}[
  x={(1.2cm,0cm)}, y={(0.55cm,0.85cm)}, z={(0cm,1.1cm)},
  pt/.style={circle, fill=black, inner sep=1.2pt},
  lbl/.style={font=\scriptsize, inner sep=1pt},
  ed/.style={line width=0.6pt},
  dotline/.style={line width=0.7pt, draw=gray!65, dash pattern=on 0.5pt off 3pt, line cap=round}
]
\coordinate (v0) at (-1,-1,-1);
\coordinate (v1) at ( 1,-1,-1);
\coordinate (v2) at ( 1, 1,-1);
\coordinate (v3) at (-1, 1,-1);
\coordinate (v4) at (-1,-1, 1);
\coordinate (v5) at ( 1,-1, 1);
\coordinate (v6) at ( 1, 1, 1);
\coordinate (v7) at (-1, 1, 1);

\draw[ed]  (v0) -- (v1) -- (v2) -- (v3) -- cycle;
\draw[ed]  (v4) -- (v5) -- (v6) -- (v7) -- cycle;
\draw[ed]  (v0) -- (v4);
\draw[ed]  (v1) -- (v5);
\draw[ed]  (v2) -- (v6);
\draw[ed]  (v3) -- (v7);

\coordinate (n8)  at ( 0,-1,-1);
\coordinate (n9)  at ( 1, 0,-1);
\coordinate (n10) at ( 0, 1,-1);
\coordinate (n11) at (-1, 0,-1);

\coordinate (n12) at ( 0,-1, 1);
\coordinate (n13) at ( 1, 0, 1);
\coordinate (n14) at ( 0, 1, 1);
\coordinate (n15) at (-1, 0, 1);

\coordinate (n16) at (-1,-1, 0);
\coordinate (n17) at ( 1,-1, 0);
\coordinate (n18) at ( 1, 1, 0);
\coordinate (n19) at (-1, 1, 0);

\coordinate (n20) at ( 0,-1, 0);
\coordinate (n21) at ( 1, 0, 0);
\coordinate (n22) at ( 0, 1, 0);
\coordinate (n23) at (-1, 0, 0);
\coordinate (n24) at ( 0, 0,-1);
\coordinate (n25) at ( 0, 0, 1);

\coordinate (n26) at (0,0,0);

\draw[dotline] (n20) -- (n22);
\draw[dotline] (n21) -- (n23);
\draw[dotline] (n24) -- (n25);

\draw[dotline] (n8)  -- (n10);
\draw[dotline] (n11) -- (n9);
\draw[dotline] (n12) -- (n14);
\draw[dotline] (n15) -- (n13);
\draw[dotline] (n8)  -- (n12);
\draw[dotline] (n16) -- (n17);
\draw[dotline] (n10) -- (n14);
\draw[dotline] (n19) -- (n18);
\draw[dotline] (n11) -- (n15);
\draw[dotline] (n16) -- (n19);
\draw[dotline] (n9)  -- (n13);
\draw[dotline] (n17) -- (n18);

\foreach \i/\v in {0/v0,1/v1,2/v2,3/v3,4/v4,5/v5,6/v6,7/v7}{
  \node[pt] at (\v) {};
  \node[lbl, anchor=west] at ($(\v)+(0.06,0.04,0)$) {\i};
}
\foreach \i/\v in {8/n8,9/n9,10/n10,11/n11,12/n12,13/n13,14/n14,15/n15,16/n16,17/n17,18/n18,19/n19}{
  \node[pt] at (\v) {};
  \node[lbl, anchor=west] at ($(\v)+(0.06,0.04,0)$) {\i};
}
\foreach \i/\v in {20/n20,21/n21,22/n22,23/n23,24/n24,25/n25}{
  \node[pt] at (\v) {};
  \node[lbl, anchor=west] at ($(\v)+(0.06,0.04,0)$) {\i};
}
\node[pt] at (n26) {};
\node[lbl, anchor=west] at ($(n26)+(0.06,0.04,0)$) {26};

\node[anchor=west] at (3,-0.1,0) {%
\scriptsize
\setlength{\tabcolsep}{3pt}
\renewcommand{\arraystretch}{1.05}
\begin{tabular}{c|cccccccc}
\toprule
child & \multicolumn{8}{c}{vertices} \\
\midrule
0 & 0 &  8 & 24 & 11 & 16 & 20 & 26 & 23 \\
1 & 8 &  1 &  9 & 24 & 20 & 17 & 21 & 26 \\
2 & 24&  9 &  2 & 10 & 26 & 21 & 18 & 22 \\
3 & 11& 24 & 10 &  3 & 23 & 26 & 22 & 19 \\
4 & 16& 20 & 26 & 23 &  4 & 12 & 25 & 15 \\
5 & 20& 17 & 21 & 26 & 12 &  5 & 13 & 25 \\
6 & 26& 21 & 18 & 22 & 25 & 13 &  6 & 14 \\
7 & 23& 26 & 22 & 19 & 15 & 25 & 14 &  7 \\
\bottomrule
\end{tabular}
};
\end{tikzpicture}

\captionof{figure}{Hex27 nodes and refinement.}
\label{fig:hex27children}
\end{minipage}

\vspace{8mm}

\begin{minipage}[t]{0.92\textwidth}
\centering
\begin{tikzpicture}[
  x={(-1.3cm,-1.5cm)},
  y={(2.8cm,-.7cm)},
  z={(0cm,2.6cm)},
  pt/.style={circle, fill=black, inner sep=1.6pt},
  ptgray/.style={circle, fill=gray!65, inner sep=1.6pt},
  lbl/.style={font=\scriptsize, inner sep=1pt},
  lblgray/.style={font=\scriptsize\color{gray!65}, inner sep=1pt},
  ed/.style={line width=0.8pt},
  dotline/.style={
    line width=0.95pt,
    draw=gray!65,
    dash pattern=on 0.6pt off 3.5pt,
    line cap=round
  }
]

\coordinate (v0) at (0,0,0);
\coordinate (v1) at (1,0,0);
\coordinate (v2) at (0,1,0);
\coordinate (v3) at (0,0,1);

\draw[ed] (v0)--(v1)--(v2)--cycle;
\draw[ed] (v0)--(v3);
\draw[ed] (v1)--(v3);
\draw[ed] (v2)--(v3);

\coordinate (n0) at (0,0,0);
\coordinate (n1) at (1,0,0);
\coordinate (n2) at (0,1,0);
\coordinate (n3) at (0,0,1);

\coordinate (n4) at (0.5,0,0);
\coordinate (n5) at (0.5,0.5,0);
\coordinate (n6) at (0,0.5,0);
\coordinate (n7) at (0,0,0.5);
\coordinate (n8) at (0.5,0,0.5);
\coordinate (n9) at (0,0.5,0.5);

\coordinate (n10) at (0.333333,0.333333,0);
\coordinate (n11) at (0.333333,0,0.333333);
\coordinate (n12) at (0.333333,0.333333,0.333333);
\coordinate (n13) at (0,0.333333,0.333333);
\coordinate (n14) at (0.25,0.25,0.25);

\draw[dotline] (n4)--(n5)--(n6)--cycle;
\draw[dotline] (n4)--(n8)--(n7)--cycle;
\draw[dotline] (n6)--(n9)--(n7)--cycle;
\draw[dotline] (n5)--(n9)--(n8)--cycle;

\draw[dotline] (n5) -- (n7);

\foreach \i/\p in {0/n0,1/n1,2/n2,3/n3,4/n4,5/n5,6/n6,7/n7,8/n8,9/n9}{
  \node[pt] at (\p) {};
  \node[lbl, anchor=west] at ($(\p)+(0.025,0.025,0)$) {\i};
}

\foreach \i/\p in {10/n10,11/n11,12/n12,13/n13}{
  \node[ptgray] at (\p) {};
  \node[lblgray, anchor=west] at ($(\p)+(0.025,0.025,0)$) {\i};
}

\foreach \i/\p in {14/n14}{
  \node[ptgray] at (\p) {};
  \node[lblgray, anchor=west] at ($(\p)+(0.095,0.025,0)$) {\i};
}

\node[anchor=west] at (3.7cm,0.75cm) {%
\scriptsize
\setlength{\tabcolsep}{3pt}
\renewcommand{\arraystretch}{1.05}
\begin{tabular}{c|cccc}
\toprule
child & \multicolumn{4}{c}{vertices} \\
\midrule
0 & 0 & 4 & 6 & 7 \\
1 & 4 & 1 & 5 & 8 \\
2 & 6 & 5 & 2 & 9 \\
3 & 7 & 8 & 9 & 3 \\
4 & 5 & 6 & 4 & 7 \\
5 & 8 & 7 & 5 & 4 \\
6 & 7 & 9 & 8 & 5 \\
7 & 9 & 5 & 7 & 6 \\
\bottomrule
\end{tabular}
};

\end{tikzpicture}

\captionof{figure}{Tet15 nodes and refinement.}
\label{fig:tet15children}
\end{minipage}

\vspace{8mm}

\begin{minipage}[t]{0.92\textwidth}
\centering
\begin{tikzpicture}[
  x={(-1.cm,-1.2cm)},
  y={(2.cm,-.5cm)},
  z={(0cm,1.6cm)},
  pt/.style={circle, fill=black, inner sep=1.2pt},
  ptgray/.style={circle, fill=gray!65, inner sep=1.2pt},
  lbl/.style={font=\scriptsize, inner sep=1pt},
  lblgray/.style={font=\scriptsize\color{gray!65}, inner sep=1pt},
  ed/.style={line width=0.6pt},
  dotline/.style={line width=0.7pt, draw=gray!65, dotted}
]


\coordinate (n0) at (0,0,-1);
\coordinate (n1) at (1,0,-1);
\coordinate (n2) at (0,1,-1);

\coordinate (n3) at (0,0, 1);
\coordinate (n4) at (1,0, 1);
\coordinate (n5) at (0,1, 1);

\coordinate (n6) at (0.5,0,-1);      
\coordinate (n7) at (0.5,0.5,-1);    
\coordinate (n8) at (0,0.5,-1);      
\coordinate (n9)  at (0.5,0, 1);     
\coordinate (n10) at (0.5,0.5, 1);   
\coordinate (n11) at (0,0.5, 1);     
\coordinate (n12) at (0,0,0);        
\coordinate (n13) at (1,0,0);        
\coordinate (n14) at (0,1,0);        

\coordinate (n15) at (0.5,0,0);            
\coordinate (n16) at (0.5,0.5,0);          
\coordinate (n17) at (0,0.5,0);            
\coordinate (n18) at (0.333333,0.333333,-1); 
\coordinate (n19) at (0.333333,0.333333, 1); 

\coordinate (n20) at (0.333333,0.333333,0);

\draw[ed] (n0)--(n1)--(n2)--cycle;          
\draw[ed] (n3)--(n4)--(n5)--cycle;          
\draw[ed] (n0)--(n3);
\draw[ed] (n1)--(n4);
\draw[ed] (n2)--(n5);

\draw[dotline] (n6)--(n7);
\draw[dotline] (n7)--(n8);
\draw[dotline] (n8)--(n6);

\draw[dotline] (n12)--(n13);
\draw[dotline] (n13)--(n14);
\draw[dotline] (n14)--(n12);

\draw[dotline] (n15)--(n16);
\draw[dotline] (n16)--(n17);
\draw[dotline] (n17)--(n15);

\draw[dotline] (n9)--(n10);
\draw[dotline] (n10)--(n11);
\draw[dotline] (n11)--(n9);

\draw[dotline] (n6)--(n9);
\draw[dotline] (n7)--(n10);
\draw[dotline] (n8)--(n11);

\foreach \i in {0,...,17}{
  \node[pt] at (n\i) {};
  \node[lbl, anchor=west] at ($(n\i)+(0.06,0.04,0)$) {\i};
}
\foreach \i in {18,19,20}{
  \node[ptgray] at (n\i) {};
  \node[lblgray, anchor=west] at ($(n\i)+(0.06,0.04,0)$) {\i};
}

\node[anchor=west] at (3.7cm,-0.5cm) {%
\scriptsize
\setlength{\tabcolsep}{3pt}
\renewcommand{\arraystretch}{1.05}
\begin{tabular}{c|cccccc}
\toprule
child & \multicolumn{6}{c}{vertices} \\
\midrule
0 & 0 &  6 &  8 & 12 & 15 & 17 \\
1 & 6 &  1 &  7 & 15 & 13 & 16 \\
2 & 8 &  7 &  2 & 17 & 16 & 14 \\
3 & 7 &  8 &  6 & 16 & 17 & 15 \\
4 & 12& 15 & 17 & 3  &  9 & 11 \\
5 & 15& 13 & 16 & 9  &  4 & 10 \\
6 & 17& 16 & 14 & 11 & 10 & 5  \\
7 & 16& 17 & 15 & 10 & 11 & 9  \\
\bottomrule
\end{tabular}
};

\end{tikzpicture}

\captionof{figure}{Wedge21 nodes and refinement.}
\label{fig:tet21children}
\end{minipage}

\end{figure}


\paragraph{Adjacency choice}
The adjacency rule used to enforce grading can be selected by the user:
(i) \emph{face/edge adjacency}, or
(ii) \emph{vertex adjacency}.
In the second case, adjacency through vertices implicitly includes adjacency through edges and faces.

In the numerical experiments presented here, face/edge adjacency is used for quadrilateral and hexahedral meshes, while vertex adjacency is adopted for triangular, tetrahedral, and wedge meshes.

\paragraph{One-level (2:1) grading constraint}
The multilevel hierarchy is required to satisfy a one-level grading constraint:
any two adjacent elements may differ in refinement level by at most one.
More precisely, if two elements share a face (and, depending on the chosen adjacency rule, also an edge or a vertex),
their refinement levels $\ell$ and $\ell'$ must satisfy
\[
|\ell - \ell'| \le 1 .
\]
This constraint prevents the occurrence of arbitrarily large level jumps between neighboring elements and ensures geometric regularity of the hierarchy under refinement.

\paragraph{Cut elements}
An element $K$ is said to be \emph{cut} if the nodal values of the level-set function $\PsiLS$ on $K$ exhibit a sign change or if at least one node satisfies $\PsiLS=0$.

\paragraph{Storage convention}
The level-set field $\PsiLS$ is stored only on the finest mesh $\TT^{\levelN}$.
During multilevel marking and refinement procedures, $\PsiLS$ is re-evaluated analytically on coarser levels as needed.

\section{Multilevel mesh construction with one-level grading} \label{sec:refinement}

\subsection{Refinement templates} \label{sec:amr}
\label{subsec:refinement_templates}

In two spatial dimensions, each refined element generates $4$ children, while in three dimensions each refinement produces $8$ children.
Refinement is defined element-wise through fixed connectivity templates, and all parent-child relations are expressed in reference coordinates of the corresponding canonical element.
Figure~\ref{fig:quad9children}-\ref{fig:tet21children} illustrate, in dashed lines, the refinement templates adopted for the element families considered in this work.

For all element families except Tet15 (namely, Quad9, Tri7, Hex27, and Wedge21), the refinement produces children that are all mutually similar, up to rigid motions (translations and rotations) and uniform scaling.
In these cases, refinement corresponds to an affine subdivision of the reference element, and each child can be obtained from a single canonical child template by composition with a rigid motion and scaling.
Consequently, a single reference element and a single parent-child connectivity template per element type suffice to describe the geometry of all children generated at any refinement level.

The tetrahedral case requires special attention.
For Tet15 elements, the refinement rule generates a finite but nontrivial set of similarity classes of child tetrahedra.
Recursive application of the refinement produces six distinct similarity families up to rigid motions and scaling, which reduce to three families if mirror equivalence is also admitted.
This behavior is well known in tetrahedral refinement schemes based on edge- or face-based subdivisions and contrasts with the uniform similarity observed for quadrilateral, triangular, hexahedral, and wedge elements \cite{Bey2000,Rivara1984,Kossaczky1994}.

Crucially, the set of tetrahedral similarity classes generated by the Tet15 refinement is finite and closed under further refinement.
No new similarity types arise beyond a fixed depth, and all child tetrahedra remain shape-regular.
All tetrahedral children, regardless of their similarity class, are mapped to the same canonical Tet15 reference element via orientation-preserving affine mappings.

The finiteness of the similarity classes implies that the refinement rule avoids generating arbitrarily distorted or degenerate tetrahedra.
As a result, the Jacobians of the geometric mappings and their inverses remain uniformly bounded away from singularity across all refinement levels \cite{Bey2000}.
This property is essential for the robustness of inverse mappings, interpolation operators, and numerical integration, and it ensures stability of the discretization under repeated adaptive refinement.


\subsection{Mark-and-closure refinement}
Given a mesh $\TT^{\ell-1}$ satisfying the one-level ($2\!:\!1$) grading constraint, refinement is performed through a
\emph{mark-and-closure} procedure designed to preserve this property.

Starting from an initial set of marked elements $\mathcal{M}\leftarrow\mathcal{M}_0$ (selected for refinement by a given marking strategy), we enlarge $\mathcal{M}$ so that the $2\!:\!1$ rule is maintained after refinement.
Since only elements in $\mathcal{M}$ are refined, refinement increases the level of each $K\in\mathcal{M}$ by one, while the level of any $K'\notin\mathcal{M}$ remains unchanged.
As a consequence, a violation of the 2:1 grading condition can only be created across an interface between a marked element $K$ and an unmarked neighbor $K'$, and only when
\[
\ell(K)-\ell(K')=1,
\]
in which case refining $K$ alone would produce a level difference of $2$.

To prevent this, the marked set is closed by iteratively adding any adjacent element $K'\notin\mathcal{M}$ for which such a configuration occurs.
The closure process terminates when no further elements need to be added.
Only after this step is the actual refinement carried out: all elements in the final marked set are refined simultaneously according to the element-type–specific refinement templates, and the corresponding father–children relations are stored.
The resulting mesh $\TT^{\ell}$ therefore satisfies the one-level grading constraint by construction.

This procedure is summarized in Algorithm~\ref{alg:closure}.

\begin{algorithm}[htbp]
\caption{Mesh refinement with one-level closure (from level $\ell-1$ to $\ell$)}
\label{alg:closure}
\begin{algorithmic}[1]
\Require Mesh $\TT^{\ell-1}$ satisfying one-level grading; initial marked set $\mathcal{M}_0$; adjacency rule $\mathcal{N}(\cdot)$.
\Ensure Refined mesh $\TT^{\ell}$ satisfying one-level grading; father/children relations updated.
\State $\mathcal{M} \gets \mathcal{M}_0$
\Repeat
  \State $\mathcal{M}_{\mathrm{new}} \gets \mathcal{M}$
  \ForAll{$K \in \mathcal{M}$}
    \ForAll{$K'\notin \mathcal{M}$ and in the neighbor of $K$}
      \If{refining $K$ would create a level jump $>1$ between $K$ and $K'$}
        \State $\mathcal{M}_{\mathrm{new}} \gets \mathcal{M}_{\mathrm{new}} \cup \{K'\}$
      \EndIf
    \EndFor
  \EndFor
  \State $\mathcal{M}_{\mathrm{old}} \gets \mathcal{M}$; \quad $\mathcal{M} \gets \mathcal{M}_{\mathrm{new}}$
\Until{$\mathcal{M} = \mathcal{M}_{\mathrm{old}}$}
\State Refine all elements in $\mathcal{M}$ according to the element-type template to obtain $\TT^{\ell}$
\State Store father/children relations for all refined elements
\end{algorithmic}
\end{algorithm}

\subsection{Local parent-child transfer operator}
For each element type, we precompute a local prolongation (coarse-to-fine) matrix on the reference element.
Let $\{\varphi_{c,i}\}$ denote the basis (test) functions of the father element, and let $\{\xi_{f,j}\}$ be the nodal points of the fine element mapped into the father reference coordinates. The local prolongation is
\[
P_{ij} = \varphi_{c,i}(\xi_{f,j}).
\]
All operators are applied element-locally; we store one $P$ per element type and reuse it throughout.

\section[Adaptive initialization of the multilevel mesh and level-set]{Adaptive initialization of the multilevel mesh and level-set $\PsiLS$} \label{sec:init}
The level-set field $\PsiLS$ is initialized analytically in the physical domain, while the adaptive multilevel mesh hierarchy is constructed by successive refinements concentrated in a neighborhood of the interface $\PsiLS=0$.

At each refinement level $\ell$, the set of elements to be refined is determined through a two-step marking strategy.
First, all cut elements are identified, namely, elements whose nodal values of $\PsiLS$ change sign or contain at least one zero node.
These elements directly approximate the interface and constitute the core of the marked set.

Second, this set is enlarged by adding one layer of adjacent elements, as defined by the selected adjacency rule.
This additional layer ensures a sufficiently thick refined region around the interface, improving the robustness of subsequent characteristic tracing and interpolation.

After this marking step, the one-level (2:1) grading constraint is enforced through the closure procedure described in Section~\ref{sec:refinement}, and all elements in the resulting marked set are refined.
The process is repeated independently at each refinement level until the maximum level is reached.

The outcome of this procedure is the initial multilevel mesh hierarchy
\[
\{\TT^{0},\TT^{1},\dots,\TT^{\levelN}\},
\]
together with nodal values of the level-set field $\PsiLS$ stored only on the finest mesh $\TT^{\levelN}$.
No level-set data are retained on coarser levels. For details, see Algorithm~\ref{alg:init}.

\begin{algorithm}[htb]
\caption{Adaptive initialization of $\PsiLS$ and construction of the multilevel hierarchy $\{\TT^0,\dots,\TT^{\levelN}\}$}
\label{alg:init}
\begin{algorithmic}[1]
\Require Uniform mesh $\TT^{0}$; maximum level $\levelN$; analytic $\PsiLS({x})$; adjacency rule.
\Ensure Multilevel hierarchy $\{\TT^{0},\TT^{1},\dots,\TT^{\levelN}\}$ with father/children relations; nodal values $\PsiLS$ stored on $\TT^{\levelN}$.
\For{$\ell = 1$ to $\levelN$}
  \State Re-evaluate $\PsiLS$ analytically at nodal coordinates of the current mesh $\TT^{\ell-1}$
  \State Mark elements $K\in \TT^{\ell-1}$ such that (i) $\PsiLS$ changes sign across nodes of $K$, or (ii) at least one node has $\PsiLS=0$
  \State Close the marked set by the neighbor-marking fixed-point procedure (one-level constraint) using Algorithm~\ref{alg:closure}
  \State Refine marked elements to obtain $\TT^{\ell}$ and store father/children relations
\EndFor
\State On $\TT^{\levelN}$, evaluate and store $\PsiLS$ at all finest-level nodes
\end{algorithmic}
\end{algorithm}

\section{Point location on a multilevel mesh} \label{sec:pointlocation}
Given a physical point ${x}$, we determine whether it lies inside the computational domain, and if so, we compute the containing element and reference coordinates ${\xi}$ such that $F_K({\xi})={x}$, where $F_K : \hat K \to K$ denotes the mapping from the reference element $\hat K$ to the physical element $K$.

\subsection{Coarse Cartesian indexing for level-0 elements}
\label{subsec:loc-coarse}
We build a Cartesian grid of boxes covering the domain with characteristic size $\approx 0.1$ times the minimum level-0 element size.
For each box, we compute an ordered list of candidate overlapping level-0 elements using bounding-box overlap and a barycenter-based ordering:
first an element containing the box barycenter (if any), followed by candidates sorted by increasing distance between the barycenter and element nodes.
If the barycenter belongs to an element, we also precompute and store the inverse mapping of the box center for that element.

To reduce the cost of repeatedly evaluating ${x}=F_K({\xi})$ and its Jacobian during the Newton iteration, we employ a polynomial representation of the isoparametric map that is uniform across all element
families present in the mesh.
For each level-0 element $K$, 
the mapping $F_K : \hat K \to K$ is expressed in a polynomial basis on the reference element as
\[
F_K({\xi})
= \sum_{i+j+k \le n} \mathbf{a}_{ijk}^{(K)}\, \xi_1^i \xi_2^j \xi_3^k .
\]
The coefficients $\mathbf{a}_{ijk}^{(K)}$ depend only on the physical coordinates of the nodes of $K$
and are precomputed once and stored. For two-dimensional elements (Quad9 and Tri7), the third reference coordinate is absent and the expansion reduces to the case $k=0$.

This representation provides a unified and efficient evaluation strategy for both $F_K({\xi})$ and
$\nabla F_K({\xi})$ across mixed element types, avoiding element-specific evaluation kernels inside the
Newton iteration.
The resulting reduction in computational cost is particularly beneficial in point-location and
particle-tracking routines, where multiple candidate elements may be tested.
Further details are given in \cite{capodaglio2017particle}.

\subsection{Inverse mapping and Newton criteria}
\label{subsec:inverse-map}

Given a candidate element $K$ and a physical point ${x}$, we compute the corresponding
reference coordinate ${\xi} = F_K^{-1}({x})$ by applying Newton’s method to the
nonlinear system
\[
F_K({\xi}) - {x} = {0},
\]
where $F_K : \hat K \rightarrow K$ denotes the isoparametric mapping from the reference element
$\hat K$ to the physical element $K$.

Starting from an initial guess ${\xi}^{(0)} \in \hat K$, the Newton iteration reads
\[
{\xi}^{(m+1)}
=
{\xi}^{(m)}
-
\left[
J_{F_K}\!\left({\xi}^{(m)}\right)
\right]^{-1}
\left(
F_K\!\left({\xi}^{(m)}\right) - {x}
\right).
\]

The iteration is terminated when the residual norm
$\|F_K({\xi}^{(m)}) - {x}\|$
falls below a prescribed tolerance $\epsilon$, or when a maximum number of iterations is reached.
If convergence is not achieved within the iteration cap, the candidate element is discarded
without performing any reference-domain inclusion test, and the next element in the ordered list
is examined. If the Newton iteration converges, the resulting reference coordinate
${\xi}$ is tested for inclusion in the reference element $\hat K$.
If ${\xi}$ lies within the admissible bounds of $\hat K$, the point
${x}$ is declared to be located in the element $K$.
Otherwise, the candidate element is rejected and the search proceeds with the
next element in the ordered list, restarting the inverse-mapping procedure.

If all candidate elements have been examined and none yield a convergent
Newton solution with ${\xi} \in \hat K$, the point ${x}$ is marked as lying
outside the computational domain.

\paragraph{Initial guess and acceleration strategy}
The initial Newton iterate ${\xi}^{(0)}$ is chosen according to the candidate-selection logic
described in Section~\ref{subsec:loc-coarse}.
If the element is the preferred candidate associated with the Cartesian box barycenter and a stored
inverse mapping is available, this value is used as ${\xi}^{(0)}$.
Otherwise, the reference coordinates of the closest element vertex provide a warm start.

To accelerate convergence for higher-order elements (Quad9, Tri7, Hex27, Tet15, Wedge21), the Newton iteration is initialized using a linearized inverse mapping.
Specifically, we first perform a small number of Newton iterations on the corresponding
linear element map obtained by restricting $F_K$ to its vertex nodes only (Quad4, Tri3, Hex8, Tet4, Wedge6).
The resulting reference coordinate serves as an improved initial guess for the full
higher-order Newton iteration.
This strategy significantly reduces the number of nonlinear iterations and improves robustness
for highly curved or distorted elements.

\subsection{Fast push-down using the multilevel hierarchy}
\label{subsec:pushdown}

Once a point is located in an element $K^\ell \in \TT^\ell$ with reference coordinate ${\xi}^\ell$, its location in finer levels can be obtained without any additional global search.
If $K^\ell$ is refined, the child element $K^{\ell+1}$ containing the point and the corresponding
reference coordinate ${\xi}^{\ell+1}$ are computed directly by applying the known affine
parent-to-child reference mappings.
This operation is purely local and incurs negligible cost.

This push-down mechanism is used in two distinct phases of the algorithm.

\paragraph{(i) Push-down during multilevel mesh reconstruction}
When building a new multilevel hierarchy at time $t^{n+1}$, markers advected forward in time are first
located on the coarse mesh $\TT^0_{n+1}$.
The element in which a marker is found determines which elements are marked for refinement at level~0.
After closure and refinement, the same marker is immediately pushed down to the refined level~1
using the parent-to-child maps.
This process is repeated level by level: at each refinement stage, the marker location identifies
elements to be marked (together with their adjacency-based closure), and the marker is then pushed
down to the next level.
In this context, point location drives multilevel mesh construction.

\paragraph{(ii) Push-down during point location on an existing hierarchy}
When locating a point in an already constructed multilevel mesh (e.g.\ during backward characteristic
tracing for the update of $\PsiLS$), the point is first located on the coarse mesh $\TT^0$ using the
Cartesian box index and inverse mapping.
Once a containing level-0 element is identified, the point is pushed down deterministically through
the existing hierarchy until the finest level $\TT^{\levelN}$ is reached.
In this case, no refinement or marking is performed: the push-down is used solely to locate the point efficiently on the finest mesh.

In both cases, the global search is performed only at level~0.
All finer-level locations are obtained by local affine mappings, which makes the procedure efficient,
robust, and well-suited for repeated semi-Lagrangian queries on highly adaptive multilevel meshes.

\paragraph{Marker location Algorithm summary} 
Algorithm~\ref{alg:locate} performs point location for a marker position \(x\) in two stages.
In the first stage, \(x\) is located on the coarse mesh \(\mathcal{T}^0\).
To this end, we identify the level--0 Cartesian box \(B\) containing \(x\) and retrieve its associated precomputed ordered list \(\mathcal{C}(B)\) of overlapping level--0 elements.
The elements in \(\mathcal{C}(B)\) are examined sequentially.

For each candidate element \(K \in \mathcal{C}(B)\), the inverse isoparametric mapping \(\xi = F_K^{-1}(x)\) is computed using Newton’s method.
The initial guess \(\xi^{(0)}\) is chosen according to the following strategy:
if the element \(K\) contains the barycenter of \(B\), \(\xi^{(0)}\) is set to the precomputed inverse mapping of the box barycenter \(\xi_B\);
otherwise, a warm start is constructed from the reference coordinates of the element node closest to \(x\).
A candidate element is accepted if the Newton iteration converges and the resulting reference coordinate satisfies the inclusion criterion \(\xi \in \hat K\).
The search terminates as soon as a valid candidate is found; if no candidate is accepted, the marker is classified as \texttt{out-of-domain}.

Once a containing element is identified at level~0, the second stage of Algorithm~\ref{alg:locate} propagates the marker through the refinement hierarchy.
The reference coordinate is mapped from parent to child elements using the known affine parent-to-child transformations, thereby locating the marker on finer levels without performing any additional global searches.

\begin{algorithm}[htb]
\caption{Locate a point ${x}$ in the hierarchy $\{\TT^0,\dots,\TT^{\levelN}\}$}
\label{alg:locate}
\begin{algorithmic}[1]
\Require Point ${x}$ in physical coordinates; level-0 box index; hierarchy $\{\TT^0,\dots,\TT^{\levelN}\}$ with father/children relations; tolerance $\epsilon$ and max iters.
\Ensure Either (element $K^{\levelN}$, reference coordinate ${\xi}^{\levelN}$) or \texttt{out-of-domain}.
\State Identify Cartesian box $B$ containing ${x}$
\State Retrieve the ordered candidate list $\mathcal{C}(B) = (K_1,K_2,\dots,K_m)$ of level-0 elements
\State $\texttt{found} \gets \texttt{false}$
\For{$i=1$ to $m$ \textbf{while} $\texttt{found}=\texttt{false}$}
  \State $K \gets K_i$
  \State Choose Newton initial guess:
  \If{$K$ is the barycenter-containing element of $B$ and a stored inverse map ${\xi}_B$ exists}
    \State ${\xi}^{(0)} \gets {\xi}_B$
  \Else
    \State ${\xi}^{(0)} \gets$ reference coordinate of the closest node (warm start)
  \EndIf
  \State Run Newton to solve $F_K({\xi})={x}$ with tolerance $\epsilon$ and max iters
  \If{Newton converged \textbf{and} ${\xi} \in \hat K$}
    \State $\texttt{found} \gets \texttt{true}$
    \State $(K^0,{\xi}^0)\gets(K,{\xi})$
  \EndIf
\EndFor
\If{$\texttt{found}=\texttt{false}$}
  \State \Return \texttt{out-of-domain}
\EndIf
\For{$\ell=0$ to $\levelN-1$}
  \State Using $(K^\ell,{\xi}^\ell)$ and the parent-to-child affine maps, determine $(K^{\ell+1},{\xi}^{\ell+1})$
\EndFor
\State \Return $(K^{\levelN},{\xi}^{\levelN})$
\end{algorithmic}
\end{algorithm}

\section{Reinitialization}
\label{sec:reinit}

\subsection{Motivation}
Although any Lipschitz continuous function $\PsiLS$ such that $\PsiLS=0$ on the interface can be used as a level-set function, it is standard to choose one with a controlled gradient magnitude near the interface. A common choice is the signed distance to the interface. In this work, we instead use a mollified signed-distance function
\begin{equation}
\label{mollifier}
S(d) = \begin{cases}
-1 & d \leq -\varepsilon \\
-1 + 2 (C_0 + d(C_1 + C_3 d^2)) & d \in (-\varepsilon, \varepsilon) \\
1 & d \geq \varepsilon
\end{cases}
\,,
\end{equation}
where $d$ is the signed distance, $\varepsilon$ is the half-bandwidth, $C_0 = \frac{1}{2}$, $C_1 = \frac{3}{4 \varepsilon}$ and $C_3 = -\frac{1}{4\varepsilon^3}$.

The reason for enforcing this controlled behavior of $\PsiLS$ is twofold. First, a bounded gradient magnitude that does not go to zero near the interface fixes the scale of $\PsiLS$ and keeps geometric quantities, namely the normal and the curvature,
\[
\widehat{n}=\frac{\nabla\PsiLS}{\|\nabla\PsiLS\|}, \qquad 
\kappa=\nabla\!\cdot\!\left(\frac{\nabla\PsiLS}{\|\nabla\PsiLS\|}\right),
\]
well-conditioned.
Second, the field advection is more robust. For the transport equation
\[
\partial_t \PsiLS + v\!\cdot\nabla\PsiLS = 0,
\]
overly steep ($\|\nabla\PsiLS\|\gg 1$) or flat ($\|\nabla\PsiLS\|\ll 1$) profiles tend to trigger spurious oscillations or excessive numerical diffusion in classical finite-element advection and tighten local CFL constraints. The same issue appears with Lagrangian (characteristic) updates, where $\tfrac{D\PsiLS}{Dt}=0$ is enforced in practice via
\[
\PsiLS^{n+1}(x^{n+1}) \approx \PsiLS^{n}\!\big(x^{n}\big),
\]
with $x^{n}$ the back-traced footpoint defined by
\[
x^{n}
= x^{n+1}
-\int_{t-\Delta t}^{t} v(x(s),s)\,ds .
\]
Because time discretization and spatial interpolation are not exact, irregular $\PsiLS$ profiles amplify the evaluation errors at $x^{n}$.

Unfortunately, the level-set gradient magnitude on the interface is not preserved under advection, because the flow kinematics stretch or compress the level–set gradient. Let
\[
\Gamma_c := \{x : \PsiLS(x,t)=c\},
\]
and introduce the strain–rate tensor
\[
\mathbf{S}(x,t) := \dfrac{1}{2}\big(\nabla {v}(x,t) + \nabla{v}(x,t)^{\!\top}\big).
\]
The component of the velocity gradient that governs the relative motion of the neighboring $\PsiLS$ iso-surfaces is its normal projection,
\[
\left.\widehat{n}^{\top}\nabla{v}\,\widehat{n}\right|_{\Gamma_c}
= \left.\widehat{n}^{\top}\mathbf{S}\,\widehat{n}\right|_{\Gamma_c},
\]
since the skew–symmetric part does not contribute. As a consequence, the magnitude of the level–set gradient evolves along characteristics according to 
\[
\frac{D}{Dt}\,\|\nabla\PsiLS\| \;=\; -\big(\widehat{n}^{\top}\mathbf{S}\,\widehat{n}\big)\,\|\nabla\PsiLS\|\,. 
\]
Hence, whenever the normal strain rate $\widehat{n}^{\top}\mathbf{S}\,\widehat{n}$ is nonzero, $\|\nabla\PsiLS\|$ departs from unity, producing local steepening or flattening even under exact transport. This purely kinematic distortion motivates periodic reinitialization to restore the field gradient. In addition, numerical diffusion during the advection process can also occur.


\subsection{Numerical approach}
In this work, we adopt a geometric strategy: the field is reconstructed in a narrow band around the interface through closest-point projection, while a cheaper extension is applied in the far field. Related closest-point formulations include hybrid projection–fast-marching approaches in finite differences \cite{chopp2001reinit,anumolu2013reinit}, full-grid projection with kink detection for robustness \cite{henri2022reinit}, and finite-element variants for linear elements \cite{parolini2007reinit}.

Our method extends this geometric idea to higher-order finite elements (i.e., bi- and tri-quadratic). 


Let $\PsiLS_h$ denote the discrete level-set function. The discrete interface is defined as
\begin{equation}
\label{eqGammaH}
\Gamma_h := \{\,x\in\Omega : \PsiLS_h(x)=0\,\},
\end{equation}
Then, we define the interface region as the union of the cut cells $\Kcut$ and their immediate neighbors $\mathcal{K}_{\mathrm{nb}}$:
\[
\Kcut := \{\, K_h \in \mathcal{T}_h \;:\; K_h \cap \Gamma_h \neq \emptyset \,\},\]
\[
\mathcal{K}_{\mathrm{nb}} := \{\, K_h \in \mathcal{T}_h \;:\; \exists\, K' \in \Kcut \text{ with } K_h \sim K' \,\},
\]
where $K_h \sim K'$ means that $K_h$ and $K'$ share an edge (2D) or a face (3D). Then the interface zone $\Omega_{int}$ and the far-field zone $\Omega_{ext}$ are defined as follows
\[
\Omega_{\mathrm{int}} := \bigcup_{K \in \Kcut \cup \mathcal{K}_{\mathrm{nb}}} K,
\qquad
\Omega_{\mathrm{ext}} := \Omega \setminus \Omega_{\mathrm{int}}.
\]

In both regions, the signed distance to the interface is computed, and the level-set field is reconstructed using the mollifier \eqref{mollifier}.

\paragraph{Interface region}
\label{sec:cutRegion}
Regarding the interface region, for each grid node, we compute the closest point on the surface by projecting it onto the zero level-set $\PsiLS = 0$. 

Specifically, we project the node using a Newton iteration. Let $x^{(0)}\in\mathbb{R}^d$ be a node in the interface region in physical coordinates. Its projection onto the zero level-set is obtained through the following update:
\begin{equation}
\label{update_proj}
    x^{(k+1)} = x^{(k)} 
    - \frac{\PsiLS(x^{(k)})}{\|\nabla \PsiLS(x^{(k)})\|^2} \, \nabla \PsiLS(x^{(k)})\,.
\end{equation} 

The projection stops when $|\PsiLS(x^k)| \le \texttt{tol}$. Ill-conditioned gradients may occur in complex two-phase flows~\cite{henri2022reinit} and can lead to either a failure to project the node or to inaccurate results. Even though we project only the interface-field nodes, which are very close to the interface and therefore less likely to suffer from poor conditioning, in the case of slow or failed convergence of the projection method, a fallback value from the far-field reconstruction is safely adopted. 
The point projection procedure is summarized in Algorithm \ref{alg::point_projection}.

\begin{algorithm}[htb]
\caption{Project nodes in interface region $\Omega_{\mathrm{int}}$ on the interface $\PsiLS = 0$}
\begin{algorithmic}[1]
\Require Hierarchy $\{\TT^0,\dots,\TT^{\levelN}\}$; finite element interpolant $\PsiLS_h$ (and $\nabla\PsiLS_h$); tolerance \texttt{tol}; max iters.
\Ensure The projection of each interface region node on the zero level-set $\PsiLS_h = 0$.
\State Collect the interface region nodes index set $\mathcal{I}$ and the corresponding coordinates $\{x_i\}_{i\in\mathcal{I}}$
\State Initialize the active index set
$\mathcal{I}_{\mathrm{active}} \gets \mathcal{I}$
\While{$it < \text{max iters}$ \textbf{and} $\mathcal{I}_{\mathrm{active}} \neq \varnothing$}
  \State Locate all points $\{x_i\}_{i\in\mathcal{I}_{\mathrm{active}}}$ in the mesh hierarchy 
    \State Evaluate the finite-element interpolant and its gradient at the located points:
    $\hspace{1cm}\{\PsiLS_h(x_i)\}_{i\in\mathcal{I}_{\mathrm{active}}}$ and $\{\nabla\PsiLS_h(x_i)\}_{i\in\mathcal{I}_{\mathrm{active}}}$
    \State $\mathcal{I}_{\mathrm{active}} \gets \{\, i\in\mathcal{I}_{\mathrm{active}} \;:\; |\PsiLS_h(x_i)| > \texttt{tol} \,\}$
    \State Update $\{x_i\}_{i\in\mathcal{I}_{\mathrm{active}}}$ using \eqref{update_proj}
\EndWhile    
\Return $\{x_i\}_{i\in\mathcal{I} \setminus \mathcal{I}_{\mathrm{active}}}$
\end{algorithmic}
\label{alg::point_projection}
\end{algorithm}

A significant portion of the numerical error introduced by the reinitialization, namely the interface shift, originates in this procedure, as it is applied in the interface region. The main source of this error is the finite-element accuracy with which the reinitialized field is represented. Indeed, at the continuous level, the signed-distance function shares the same zero level-set as the original level-set field. At the discrete level, however, the signed distance is reconstructed only at the mesh nodes and then subsequently interpolated using the finite-element basis. As a consequence, the original discrete interface defined in \eqref{eqGammaH} does not generally coincide with the reconstructed interface
\[
\widetilde{\Gamma}_h := \{\,x\in\Omega : S_h(x)=0\,\},
\]
where $S_h$ denotes the finite-element interpolant of the nodal signed-distance (or mollified signed-distance, defined in \eqref{mollifier}) values. Nevertheless, the mismatch between $\Gamma_h$ and $\widetilde{\Gamma}_h$ decreases under mesh refinement, with a convergence rate dictated by the finite-element approximation. 

\paragraph{Far-field region}
Regarding the far-field region, instead of applying a Fast Marching Method (FMM) to the remaining nodes, we perform a direct nearest-neighbor search using a kd-tree algorithm \cite{blanco2014nanoflann}. The search is carried out on a set of interface markers, whose placement is described in Section~\ref{sec:initMarkers}. This approach exhibits a similar computational scaling to the FMM. In fact, given $N_n$ cell nodes and $N_m$ markers on the surface, the computational cost of building the tree scales as $\mathcal{O}(N_m \log N_m)$, while the total cost of the searches scales as $\mathcal{O}(N_n \log N_m)$\cite{shakoor2015kdtree}. In comparison, the overall complexity of the FMM scales as $\mathcal{O}(N_n \log N_n)$. Moreover, we point out that, given our refinement strategy, the number of points $N_n$ in the far-field region for which the search must be performed is minimal.

\subsection{Marker placement}
\label{sec:initMarkers}
For each cut element $K\in\Kcut$, we construct a set of physical markers $X_{\mathrm{local}}$ placed on the zero level-set starting from its intersections with the cell's edges. 

In two dimensions, for a simple cut (characterized by two intersection points), we place the markers along the segment connecting them.  In three dimensions, we consider the cut to be simple when the number of edge intersections $N_{\mathrm{int}} \geq 3 \, \wedge N_{\mathrm{int}} \leq 5$. In this case we first compute the barycenter of those and then form a fan triangulation by creating $N_{\mathrm{int}}$ triangles sharing the barycenter as a common vertex. Markers are subsequently distributed over these triangles. For hexahedral cells, configurations with $N_{\mathrm{int}} = 6$ may still be simple. However, to avoid the geometrical check they are automatically classified as complex, without loss of generality.

Whenever the cut is no longer simple (e.g., more than two intersections in 2D or more than five intersections in 3D), we instead generate a localized marker network concentrated around the interface within the cell, by connecting the intersection points. This approach is not the most efficient in terms of the number of markers; however, while avoiding the explicit handling of complex configurations (especially in 3D), it does not compromise the accuracy of the advection and refinement routines. Moreover, the additional markers introduce a negligible overhead during advection and point-location, since such complex-cut scenarios occur only rarely.

\begin{figure}[h]
    \centering
    \includegraphics[width=0.5\linewidth]{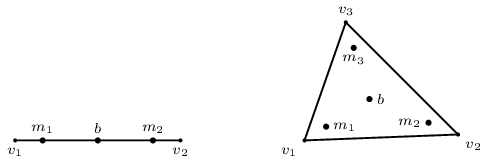}
    \caption{Advection markers placement for a 2D segment (left) and a 3D triangle(right).}
    \label{fig:advection_markers}
\end{figure}

If markers are required only for advection, we place 3 markers on each segment and 4 markers on each triangle that forms the fan triangulation: one at the barycenter $b$, and the others at positions obtained by shifting each vertex $v_i$ toward $b$ by one third of the vertex-to-barycenter distance:
\[
m_i = v_i + \frac{1}{3}(b-v_i)
    = \frac{2}{3}v_i + \frac{1}{3}b.
\]
For segments, $i=1,2$ and for triangles, $i=1,2,3$. A representation is provided in Figure~\ref{fig:advection_markers}.

Instead, when markers are used to reinitialize the level-set field, we employ a slightly higher marker density. For a prescribed linear density, markers are placed uniformly along each two-dimensional segment. For three-dimensional triangles, we take the longest edge as the base and fill the triangle with rows of markers parallel to the base. Both the marker spacing within each row and the spacing between rows are set using the same linear density. The marker density ranges from $\frac{10}{h}$ in 2D to $\frac{5}{h}$ in 3D, where $h$ is the cell size.

Finally, the resulting markers can eventually be projected onto the level-set interface using the same procedure described in Algorithm \ref{alg::point_projection}.

\subsection{Adaptive frequency} \label{sec:adapt_frequency}
To avoid unnecessary, costly reinitializations that may introduce spurious interface displacement, we monitor how far the local gradient magnitude deviates from its target value near the interface. Given the mollifier width $\varepsilon$ and expression \eqref{mollifier}, we set the target gradient magnitude to
\[
t_\varepsilon := \frac{3}{2\varepsilon}\,.
\]

We define a symmetric local deviation measure
\[
D(x) \;:=\; \left|\log\!\left(\frac{\|\nabla\PsiLS(x)\|}{t_\varepsilon}\right)\right|,
\qquad x\in \Gamma,
\]
which vanishes when $\|\nabla\PsiLS(x)\|=t_\varepsilon$ and increases for both under- and over-shoots.
Let $X=\{m_1,\dots,m_N\}$ be the set of $N_m$ marker points used for the advection step (sampling locations on the interface).
Collecting the sampled deviations into the vector
\[
\mathbf{D} \;:=\; \bigl(D(m_1),\dots,D(m_{N_m})\bigr) \in \mathbb{R}^{N_m\times d},
\]
we define the global deviation indicator as the normalized $\ell^1$-norm
\[
\overline{D}
\;:=\; \frac{1}{N_m}\,\|\mathbf{D}\|_1
\;=\; \frac{1}{N_m}\sum_{i=1}^{N_m} D(m_i).
\]
In this work we use the $\ell^1$ norm, but other $\ell^p$ norms may be adopted depending on the application.

A reinitialization step is activated if $\overline{D}$ is higher of a given tolerance $\tau_{\mathrm{reinit}}$, and skipped otherwise.
This criterion activates reinitialization only when significant deviations of $\|\nabla\PsiLS\|$ from $t_\varepsilon$ occur, thereby reducing both computational cost and the risk of unnecessary interface motion.

%

\section{One adaptive semi-Lagrangian time step}
\label{sec:timestep}

We describe one complete time step $$\text{from } (\TTn,\PsiLS^n)\rightarrow (\TTnp,\PsiLS^{n+1}),$$ where both input and output include the full multilevel mesh hierarchy.
Time is discretized with a fixed time step $\dt>0$, and we denote $t^n = n\dt$.

The evolution over one time step is driven by a prescribed velocity field \(v(x,t)\),
assumed to be Lipschitz continuous in space uniformly in time, i.e.,
\[
v \in L^\infty(0,T; W^{1,\infty}(\Omega)^d),
\]
and analytically extended outside the computational domain when needed.
Under this assumption, characteristic curves are well defined and unique on each time interval
\([t^n,t^{n+1}]\), allowing both forward and backward characteristic tracing.

In the present work, the velocity field is prescribed analytically.
However, the algorithm is designed to be coupled with a multiphase flow solver, in which case the velocity is typically defined only at mesh nodes.
From a conceptual standpoint, the overall procedure remains unchanged.
The only difference arises in the evaluation of the velocity during forward or backward characteristic tracing, where intermediate Runge--Kutta stages require velocity values at points located inside mesh elements rather than at nodes.
In this setting, the same point-location and inverse-mapping strategy introduced above can be applied to evaluate the velocity at the required off-nodal positions.

Each time step is decomposed into two logically distinct phases:
a forward advection phase used exclusively to reconstruct an adaptive multilevel mesh at time $t^{n+1}$,
and a backward semi-Lagrangian phase used to update the level-set field on the new mesh. Finally, we test whether the level-set $\PsiLS$ needs to be reinitialized, and if so, we reinitialize it.

\subsection{Forward marker advection and hierarchy reconstruction}
\label{subsec:forward}

At time $t^n$ on the finest mesh $\TT^{\levelN}_n$, we identify all cut elements, i.e., elements whose nodal values of $\PsiLS^n$
change sign or contain at least one zero node.
From these cut elements, we construct a set of physical marker points placed on the zero level-set, following the procedure outlined in Section~\ref{sec:initMarkers}.

All markers are advected forward in time over one step $\dt$ using a fourth-order Runge-Kutta scheme applied to the velocity field $v(x,t)$.
The resulting marker positions at time $t^{n+1}$ are then used to reconstruct a new multilevel hierarchy $\TTnp$ starting from a freshly initialized level-$0$ mesh.

Specifically, the advected markers are first located on $\TT^0_{n+1}$ using the level-$0$ Cartesian box index and inverse-mapping procedure described in section \ref{subsec:loc-coarse}.
The elements in which the markers are found determine the initial marked set for refinement.
At each refinement level, the marked set is enlarged by the one-level (2:1) closure procedure, and all marked elements are refined simultaneously according to the element-type refinement templates.
After refinement, the markers are pushed down to the next level using the known parent-to-child affine reference mappings, without any additional global search.

This process is repeated level by level until the finest level $\TT^{\levelN}_{n+1}$ is constructed.
The outcome of the forward phase is a new multilevel mesh hierarchy at time $t^{n+1}$, adapted to the advected interface geometry,
together with updated father-children relations.
No level-set values are transported or interpolated during this phase; forward characteristic tracing is used solely for geometric purposes.

\begin{figure}[htbp]
    \centering
    \includegraphics[width=0.32\linewidth]{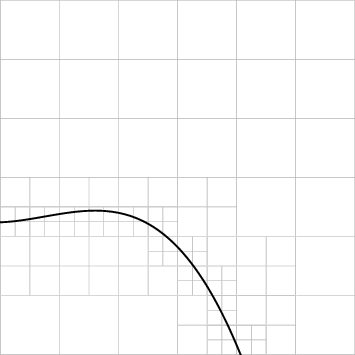}
    \includegraphics[width=0.32\linewidth]{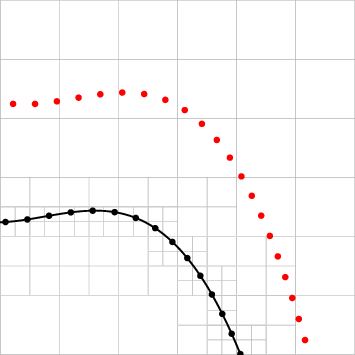}
    \includegraphics[width=0.32\linewidth]{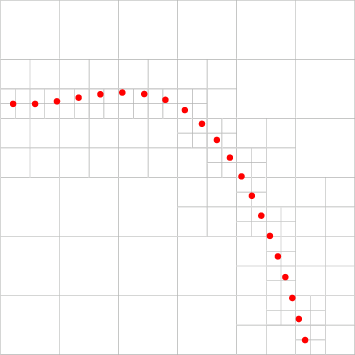}
    \caption{Schematic representation of the adaptive grid update during one advection step:
    (a) initial mesh and level-set field at time $t^n$;
    (b) marker points placed on the zero level-set and transported forward in time;
    (c) adapted mesh $\mathcal{T}^{n+1}$ refined around the advected marker positions.}
    \label{fig:marker_advection}
\end{figure}

A schematic of the forward step is given in Figure~\ref{fig:marker_advection}.

\subsection{Backward semi-Lagrangian update of the level-set field}
\label{subsec:backward}

Once the new multilevel hierarchy $\TTnp$ has been constructed, the level-set field is updated by a backward semi-Lagrangian procedure.
The arrival points are chosen as the physical coordinates of all nodes of the finest mesh $\TT^{\levelN}_{n+1}$.

For each such node $x^{n+1}$ at time $t^{n+1}=t^n+\dt$, a departure point $x^{n}$ is computed by integrating the characteristic backward in time over one step $\dt$.
Backward tracing is performed using a fourth-order Runge-Kutta scheme with time step $-\dt$, applied to the same velocity field $v(x,t)$.
This yields a high-order accurate approximation of the characteristic departure point and allows time steps that are not constrained by the mesh size in the usual CFL sense, although accuracy still depends on $\dt$.

The departure point $x^{n}$ is then located in the old multilevel hierarchy $\TTn$ using the point-location algorithm described in Algorithm~\ref{alg:locate},
without triggering any refinement.
If a containing element is found, the level-set value at the departure point is evaluated by finite-element interpolation using the nodal values stored on the deeper mesh $\TT^{\levelN}_n$. 
If the departure point lies outside the computational domain, a prescribed exterior value is assigned; in this work, we set $\PsiLS=-1$.

The semi-Lagrangian update is defined by enforcing constancy of the level-set field along characteristics,
\[
\PsiLS^{n+1}(x^{n+1}) = \PsiLS^{n}(x^{n}),
\]
for all nodes $x^{n+1}$ of $\TT^{\levelN}_{n+1}$.
The updated level–set function $\PsiLS^{n+1}$ is tested against the reinitialization algorithm described in Section~\ref{sec:reinit}; if it is excessively distorted, it is reinitialized.

After completion of the update, the level-set field $\PsiLS^{n+1}$ is stored only at the nodes of the finest mesh $\TT^{\levelN}_{n+1}$.
Coarser levels do not carry independent representations of $\PsiLS$.
The multilevel hierarchy is thus used for geometric adaptation, refinement relations, and point location,
while the level-set degrees of freedom are maintained exclusively on the active finest space.

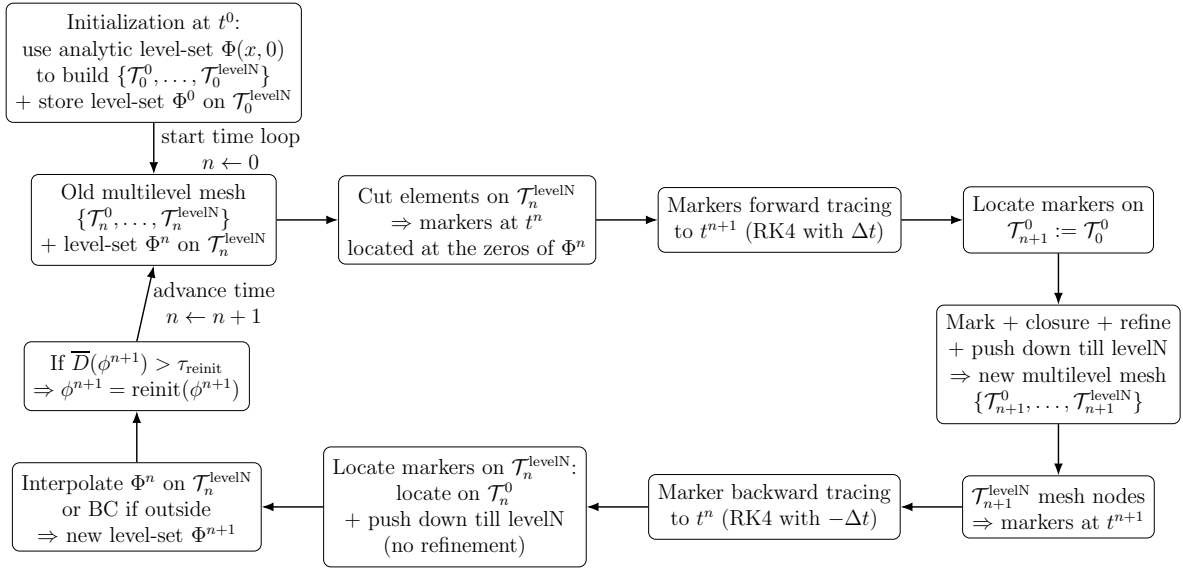
\begin{figure}[htb]
\centering
\scalebox{0.68}{%
\begin{tikzpicture}[
  node distance=10mm and 12mm,
  box/.style={rectangle, rounded corners, draw, align=center, inner sep=5pt},
  arr/.style={-Latex, thick}
]

\node[box] (init) {Initialization at $t^0$:\\
use analytic level-set $\PsiLS(x,0)$\\
to build $\{\TT^0_0,\dots,\TT^{\levelN}_0\}$\\
+ store level-set $\PsiLS^0$ on $\TT^{\levelN}_0$};

\node[box, below=of init] (old) {Old multilevel mesh\\
$\{\TT^0_n,\dots,\TT^{\levelN}_n\}$\\
+ level-set $\PsiLS^n$ on $\TT^{\levelN}_n$};

\node[box, right=of old] (markers) {Cut elements on $\TT^{\levelN}_n$\\
$\Rightarrow$ markers at $t^{n}$\\located at the zeros of $\PsiLS^n$};
\node[box, right=of markers] (fwd) {Markers forward tracing\\
to $t^{n+1}$ (RK4 with $\dt$)};

\node[box, right =of fwd] (locf) {Locate markers on\\
$\TT^0_{n+1}:= \TT^0_0$};

\node[box, below=of locf] (newmesh) {Mark + closure + refine \\ 
+ push down till levelN\\
$\Rightarrow$ new multilevel mesh\\
$\{\TT^0_{n+1},\dots,\TT^{\levelN}_{n+1}\}$
};

\node[box, below=of newmesh] (nodes) {$\TT^{\levelN}_{n+1}$ mesh nodes\\
$\Rightarrow$ markers at $t^{n+1}$};
\node[box, left=of nodes] (bwd) {Marker backward tracing\\
to $t^n$ (RK4 with $-\dt$)};
\node[box, left=of bwd] (loco) {Locate markers on $\TT^{\levelN}_{n}$:\\
locate on $\TT^{0}_{n}$ \\+ push down till $\levelN$ \\ (no refinement)};
\node[box, left=of loco] (interp) {Interpolate $\PsiLS^n$ on 
$\TT^{\levelN}_{n}$\\
or BC if outside\\
$\Rightarrow$ new level-set $\PsiLS^{n+1}$};
\node[box, above=of interp] (reinit) { 
If $\overline{D}(\phi^{n+1}) >\tau_{\text{reinit}}$  \\ $\Rightarrow \phi^{n+1}= \text{reinit}(\phi^{n+1})$ };

\draw[arr] (init.south) -- node[right,align=center]{start time loop\\ $n \leftarrow 0$} (old.north);

\draw[arr] (old.east) -- (markers.west);
\draw[arr] (markers.east) -- (fwd.west);

\draw[arr] (fwd.east) -- (locf.west);

\draw[arr] (locf.south) -- (newmesh.north);

\draw[arr] (newmesh.south) -- (nodes.north);
\draw[arr] (nodes.west) -- (bwd.east);
\draw[arr] (bwd.west) -- (loco.east);
\draw[arr] (loco.west) -- (interp.east);
\draw[arr] (interp.north) -- (reinit.south);

\draw[arr] (reinit.north) --
node[right,align=center]{advance time\\ $n\leftarrow n+1$} (old.south);

\end{tikzpicture}%
}
\caption{Overall workflow of the adaptive semi-Lagrangian scheme.}
\label{fig:recap}
\end{figure}

The complete forward-backward procedure for one adaptive semi-Lagrangian time step is summarized in Algorithm~\ref{alg:timestep} and in Figure \ref{fig:recap}.

\begin{algorithm}[t!]
\caption{One adaptive semi-Lagrangian time step (multilevel input/output)}
\label{alg:timestep}
\begin{algorithmic}[1]
\Require Multilevel hierarchy $\TTn$ with father/children relations; nodal $\PsiLS^n$ stored on $\TT^{\levelN}_n$; timestep $\dt$; analytic velocity field $v(x,t)$ (extended outside the domain); adjacency rule.
\Ensure Multilevel hierarchy $\TTnp$ with father/children relations; nodal $\PsiLS^{n+1}$ stored on $\TT^{\levelN}_{n+1}$.
\Statex \textbf{(A) Build the new hierarchy using forward-advected cut-element markers}
\State Identify cut elements $\Kcut \subset \TT^{\levelN}_n$ (sign change or zero node)
\State Build the marker set, $X^n$, of physical marker points placed on the zero level-set for all $K\in\Kcut$
\State Advect each marker in $X^n$ forward by RK4 over $\dt$ to obtain $X^{n+1}$
\State Initialize a fresh level-0 mesh $\TT^{0}_{n+1}$ (same coarse topology as $\TT^{0}_n$)
\State Locate markers in $X^{n+1}$ on $\TT^{0}_{n+1}$ using the level-0 box index (Algorithm~\ref{alg:locate}, level-0 stage)
\For{$\ell=1$ to $\levelN$}
  \State Mark elements of $\TT^{\ell-1}_{n+1}$ using located markers and apply one-level closure
  \State Refine to obtain $\TT^{\ell}_{n+1}$ and store father/children relations
  \State Push markers down to level $\ell$ using parent-to-child affine maps (no global search)
\EndFor
\Statex \textbf{(B) Update $\PsiLS$ on the new finest level using backward tracing in the old hierarchy}
\ForAll{nodes $x^{n+1}$ of $\TT^{\levelN}_{n+1}$}
  \State Trace $x^{n+1}$ backward by $\dt$ to get departure point $x^{n}$
  \State Locate $x^{n}$ in the old hierarchy $\TTn$ via Algorithm~\ref{alg:locate} (no refinement)
  \If{departure point is \texttt{out-of-domain}}
    \State $\PsiLS^{n+1}(x^{n+1}) \gets -1$ 
  \Else
    \State Interpolate $\PsiLS^n$ at $x^{n}$ in the containing element
    \State $\PsiLS^{n+1}(x^{n+1}) \gets \PsiLS^n(x^{n})$
  \EndIf
\EndFor
\State Compute the distortion indicator $\overline{D}$ for $\PsiLS^{n+1}$ (Section~\ref{sec:adapt_frequency})
\If{$ \overline{D}(\PsiLS^{n+1}) > \tau_{\mathrm{reinit}}$}
  \State $\PsiLS^{n+1} \gets \mathrm{Reinit}(\PsiLS^{n+1})$ 
\EndIf

\end{algorithmic}
\end{algorithm}

\section{Numerical results and validation methodology}
\label{sec:numres}

\subsection{Overview}
The numerical validation of the proposed level-set framework is carried out through a series of canonical kinematic benchmarks. These benchmarks are based on analytic, divergence-free velocity fields and are widely adopted in the literature for assessing the accuracy and robustness of level-set formulations. They provide a controlled environment for quantitatively evaluating two fundamental aspects of interface-capturing methods: 
(i) \emph{mass (or volume) conservation}, and 
(ii) \emph{geometric accuracy} during severe interface deformation and \eugenio{near-topological changes}.

In all tests, the velocity field reverses exactly in time, so that the interface should return to its initial configuration at the final time $t_1$. This property allows a direct quantitative comparison between the initial and final states, yielding exact reference data for the evaluation of numerical errors.

\eugenio{Because the prescribed velocity fields considered in these benchmarks are smooth and divergence-free, the corresponding continuous evolution preserves interface topology. For this reason, deformation benchmarks of this type are commonly used in the level-set literature to assess robustness under severe stretching, thin filament formation, and near-topological changes, rather than to model exact physical breakup or coalescence phenomena \cite{Rider1998, Enright2002, Wang2009, Cervone2009, Ramanuj2019}. In practice, genuine interface merging and breakup are typically investigated in fully coupled multiphase flow simulations, where the interface evolution is driven by the Navier-Stokes dynamics together with surface tension and other physical effects \cite{Wang2009, Xia2023, aulisa2025mixed}.}

\begin{figure}[!ht]
    \centering
    \begin{subfigure}[b]{0.27\linewidth}
        \centering
        \includegraphics[width=\linewidth]{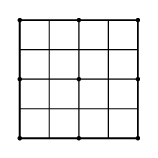}
        \caption{}
        \label{fig:box2d}
    \end{subfigure}
    \hfill
    \begin{subfigure}[b]{0.27\linewidth}
        \centering
        \includegraphics[width=\linewidth]{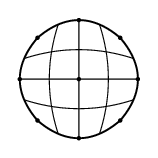}
        \caption{}
        \label{fig:ball2d}
    \end{subfigure}
    \hfill
    \begin{subfigure}[b]{0.27\linewidth}
        \centering
        \includegraphics[width=\linewidth]{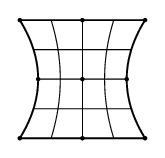}
        \caption{}
        \label{fig:funnel2d}
    \end{subfigure}
    
    \vspace{3mm}
    
    \begin{subfigure}[b]{0.27\linewidth}
        \centering
        \includegraphics[width=\linewidth]{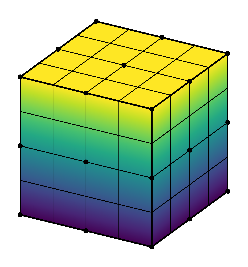}
        \caption{}
        \label{fig:box3d}
    \end{subfigure}
    \hfill
    \begin{subfigure}[b]{0.27\linewidth}
        \centering
        \includegraphics[width=\linewidth]{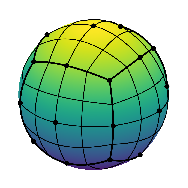}
        \caption{}
        \label{fig:ball3d}
    \end{subfigure}
    \hfill
    \begin{subfigure}[b]{0.27\linewidth}
        \centering
        \includegraphics[width=\linewidth]{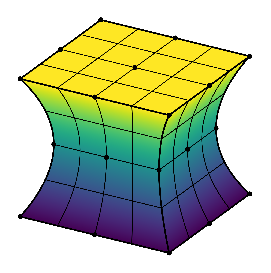}
        \caption{}
        \label{fig:funnel3d}
    \end{subfigure}
    
    \caption{Schematic view of the finite element reference meshes, including two levels of uniform refinement indicated by thin lines. Top row (2D): (a) Quad9 box element, (b) Quad9 circular element, (c) Quad9 funnel-shaped element. Bottom row (3D): (d) Hex27 box element, (e) Hex27 spherical element, (f) Hex27 funnel-shaped element. In each panel, the finite element nodes and edges are shown.}
    \label{fig:grids}
\end{figure}

In the following sections, we evaluate the above metrics on a collection of two- and three-dimensional benchmark problems, including both Cartesian and curved grids (Figure~\ref{fig:grids}), and employing all element types introduced in the previous sections.
This comprehensive set of tests allows us to assess not only the accuracy and mass conservation properties of the method, but also its overall robustness and versatility across different mesh types and problem configurations.
In Figure~\ref{fig:grids}, we display the Quad9 and Hex27 mesh configurations after two uniform refinement steps.
The Tri7 and Wedge21 meshes are obtained by subdividing each quadrilateral coarse element into two triangles and each hexahedral coarse element into two wedges, respectively.
The Tet15 meshes are fully unstructured: the box configuration is generated by partitioning the domain into five Tet15 elements, while the spherical configuration is obtained by decomposing the domain into octants.

\subsection{Error metrics and AMR-based convergence characterization}
\label{secErr}

Two complementary error measures are employed to assess the accuracy of the level-set evolution: the \emph{mass-conservation error} and the \emph{geometric error}. 

\paragraph{Mass conservation error}
The mass conservation error quantifies deviations in the total area (or volume) enclosed by the zero level-set between the initial ($t_0$) and final ($t_1$) configurations. Following the conceptual framework presented in~\cite{Rider1998,Aulisa2003Mixed}, the discrete formulation used in the present work reads:
\begin{equation}
E_{m}(t_1) = 
\dfrac{\left|\sum_{i=1}^{N_{el}} A_i C_i(t_1) - A_i C_i(t_0)\right|}
{\sum_{i=1}^{N_{el}} A_i C_i(t_0)}\,.
\label{eq_mass_error}
\end{equation}
Here, $C_i(t)$ is the phase indicator in cell $i$, i.e., the fraction of the cell occupied by the dispersed phase, and $A_i$ is the corresponding element area (in 2D) or volume (in 3D).

\paragraph{Geometric error}
The geometric accuracy is assessed by directly comparing, for each cell, the phase indicator field in the initial and final configurations. The discrete metric adopted here is defined as
\begin{equation}
E_{g}(t_1) =
\sum_{i=1}^{N_{el}} A_i \left| C_i(t_1) - C_i(t_0) \right|\,.
\label{eq_geom_error}
\end{equation}
 

The error metrics introduced above are evaluated on a sequence of successively refined meshes in order to assess the observed order of convergence for the different test cases and finite-element discretizations.
For a generic error quantity \(E\), the convergence rate between two consecutive refinement levels is defined as
\begin{equation}
p(E)
=
\frac{\log(E_i / E_{i+1})}{\log(h_i / h_{i+1})}
=
\log_2(E_i / E_{i+1}),
\label{eq_order_conv}
\end{equation}
where \(E_i\) and \(E_{i+1}\) denote the error values corresponding to mesh sizes \(h_i\) and \(h_{i+1}=h_i/2\), respectively. 

In addition, an average convergence rate is reported to reduce the influence of level-to-level oscillations,
\begin{equation}
\bar{p}(E)
=
\frac{\log\!\big(E_{\ell_{\min}}/E_{\ell_{\max}}\big)}{\log\!\big(h_{\ell_{\min}}/h_{\ell_{\max}}\big)}
=
\frac{\log_2\!\big(E_{\ell_{\min}}/E_{\ell_{\max}}\big)}{\ell_{\max}-\ell_{\min}},
\label{eq_order_conv_average}
\end{equation}
where \(\ell_{\min}\) and \(\ell_{\max}\) denote the coarsest and finest refinement levels considered.

The averaged rate provides a more reliable indicator of the asymptotic convergence behavior, particularly for error measures such as the mass error, which may exhibit non-monotone behavior due to fortuitous cancellation of local gain and loss contributions.

The same formulas also apply to the scaling of the computational time, as is shown in the next section.

\paragraph{AMR characterization, convergence and computational scaling}

In AMR computations, the numerical error depends on the local refinement pattern rather than on a uniform grid spacing. Convergence is therefore characterized by the maximum refinement level $\ell_{\max}$, which sets the minimum spacing $h_{\min}=h_{0}2^{-\ell_{\max}}$, and by the number of active leaf cells $N_{\text{cells}}$, measuring computational cost. 
Here $h_{0}$ is the characteristic length of the uniform mesh at $\ell_0$.
Accordingly, mass and geometric errors are reported as functions of $h_{\min}$ (or $\ell_{\max}$), allowing direct comparison with uniform-grid results at the same effective resolution.

In addition, from a computational viewpoint, AMR refines only a narrow band around the interface, substantially reducing the growth of degrees of freedom compared with uniform refinement. While uniform grids scale with the domain measure, $N_{\text{uniform}}\sim\mathcal{O}(h^{-d})$ in $d$ dimensions, AMR scales with the interface measure, yielding $N_{\text{AMR}}^{2D}\sim\mathcal{O}(h^{-1})$ and $N_{\text{AMR}}^{3D}\sim\mathcal{O}(h^{-2})$. This reduction leads to significantly lower execution times while preserving comparable accuracy. The scaling behavior observed in Tables~\ref{tab:2d_scaling_amr_uniform} and \ref{tab:3d_scaling_amr_uniform} confirms these theoretical trends, where execution times [$s$] and empirical scaling orders $p$ (computed via Eq.~\eqref{eq_order_conv}) are reported.

\begin{table}[h!]
\centering
\caption{Scaling comparison between AMR and uniform refinement on the Quad9 box element domain.}
\label{tab:2d_scaling_amr_uniform}
\renewcommand{\arraystretch}{1.15}
\begin{tabular}{c c c c c c c c}
\hline
$\ell_{\max}$ & $h_{\min}$ &
$N_{\text{cells}}^{\mathrm{AMR}}$ &
$t^{\mathrm{AMR}}$ &
$p^{\mathrm{AMR}}$ &
$N_{\text{cells}}^{\mathrm{uniform}}$ &
$t^{\mathrm{uniform}}$ &
$p^{\mathrm{uniform}}$ \\
\hline
6 & $2^{-6}$ & $1.76\mathrm{E}{+03}$ & 0.89 & --- & $4.10\mathrm{E}{+03}$ & 1.70 & --- \\
7 & $2^{-7}$ & $4.74\mathrm{E}{+03}$ & 2.36 & 1.40 & $1.64\mathrm{E}{+04}$ & 6.70 & 2.04 \\
8 & $2^{-8}$ & $1.22\mathrm{E}{+04}$ & 5.98 & 1.34 & $6.55\mathrm{E}{+04}$ & 28.78 & 2.04 \\
9 & $2^{-9}$ & $2.82\mathrm{E}{+04}$ & 13.88 & 1.22 & $2.62\mathrm{E}{+05}$ & 117.25 & 2.03 \\
10 & $2^{-10}$ & $6.05\mathrm{E}{+04}$ & 31.64 & 1.19 & $1.05\mathrm{E}{+06}$ & 475.21 & 2.02 \\
\hline
\end{tabular}
\end{table}

\begin{table}[h!]
\centering
\caption{Scaling comparison between AMR and uniform refinement on the cubic domain with Hex27 elements.}
\label{tab:3d_scaling_amr_uniform}
\renewcommand{\arraystretch}{1.15}
\begin{tabular}{c c c c c c c c}
\hline
$\ell_{\max}$ & $h_{\min}$ &
$N_{\text{cells}}^{\mathrm{AMR}}$ &
$t^{\mathrm{AMR}}$ &
$p^{\mathrm{AMR}}$ &
$N_{\text{cells}}^{\mathrm{uniform}}$ &
$t^{\mathrm{uniform}}$ &
$p^{\mathrm{uniform}}$ \\
\hline
4 & $2^{-4}$  & $1.44\mathrm{E}{+03}$ & 2.01   & ---  & $4.10\mathrm{E}{+03}$     & 7.05    & ---  \\
5 & $2^{-5}$  & $5.90\mathrm{E}{+03}$ & 10.89  & 2.44 & $3.28\mathrm{E}{+04}$     & 58.22   & 3.05 \\
6 & $2^{-6}$  & $2.45\mathrm{E}{+04}$ & 50.27  & 2.21 & $2.62\mathrm{E}{+05}$     & 471.60  & 3.02 \\
7 & $2^{-7}$  & $1.10\mathrm{E}{+05}$ & 228.17 & 2.18 & $2.10\mathrm{E}{+06}$     & 3787.67 & 3.01 \\
\hline
\end{tabular}
\end{table}

\subsection{Tests overview}
\label{sec:2dtests}
The two- and three-dimensional benchmarks presented in the following section are: a \emph{vortex deformation}, a \emph{rising-bubble-type} flow, and a \emph{solid-body rotation}. Each case is examined on both a standard \emph{square/cubic domain} and a \emph{circular/ball domain} in order to assess robustness against grid distortion. Simulations are carried out with maximum refinement levels in the range $\ell_{\max}=8$--$12$ for 2D and $\ell_{\max}=6$--$9$ for 3D elements. For some cases, results are compared for runs \emph{with} and \emph{without} reinitialization of the level-set field, in order to evaluate the effectiveness of the reinitialization algorithm. The accuracy is quantified using the previously defined mass and geometric errors $(E_\text{m},E_\text{g})$, while computational performance is analyzed through the scaling of total execution time with refinement level.

\subsection{Vortex deformation test}
\label{sec:vortex}

\paragraph{Test definition}
The single-vortex deformation test, introduced by LeVeque~\cite{leveque1996high}, is a stringent benchmark for level-set methods, characterized by severe interface stretching and subsequent recovery. It models a reversible kinematic deformation driven by a smooth, divergence-free velocity field.

In two dimensions, the velocity field is defined as
\begin{align}
u(x,y,t) &= \cos^2(\pi x)\,\sin(2\pi y)\,
   \cos\!\left({\pi t/T}\right),\\
v(x,y,t) &=  -\,\cos^2(\pi y)\,\sin(2\pi x)\,
   \cos\!\left({\pi t}/{T}\right),
\end{align}
with $T=8$ denoting the full advection period. The initial interface is a circle of radius $r=0.15$ centered at $(0.0,0.25)$ within $\Omega=[-0.5,0.5]^2$. At $t=T/2$, the interface is stretched into a thin filament, and at $t=T$ it should ideally recover its original configuration, allowing the computation of the mass and geometric errors.

The three-dimensional vortex deformation test generalizes the classical two-dimensional single-vortex flow to three dimensions. 
Three-dimensional vortex deformation tests of spherical interfaces, distinct from the Enright benchmark, have been previously investigated, for instance, by Aulisa et al.~\cite{aulisa2004surface,aulisa2007interface}.

We adopt the smooth, divergence-free, time-dependent velocity field used in our codebase,
\begin{align}
u(x,y,z,t) &= \cos^2(\pi x)\,\big(\sin(2\pi y)-\sin(2\pi z)\big)\,
\cos\!\left({\pi t}/{T}\right),\\[6pt]
v(x,y,z,t) &= \cos^2(\pi y)\,\big(\sin(2\pi z)-\sin(2\pi x)\big)\,
\cos\!\left({\pi t}/{T}\right),\\[6pt]
w(x,y,z,t) &= \cos^2(\pi z)\,\big(\sin(2\pi x)-\sin(2\pi y)\big)\,
\cos\!\left({\pi t}/{T}\right),
\end{align}
with T=4, defined on the periodic domain $\Omega=[-0.5,0.5]^3$.
The initial interface is a sphere of radius $r=0.15$ centered at $(0.0, 0.0, 0.25)$.
At $t=T/2$ the interface undergoes extreme filamentation, while at $t=T$ it should ideally recover its initial configuration.


\paragraph{Reinitialization study (Quad9)}
Table~\ref{tab:vortex_no_vs_reinit} reports numerical results for the 2D vortex benchmark obtained without and with reinitialization on Quad9 elements, with $\tau_{\mathrm{reinit}} = 0.25$ in the latter case. This test is chosen for the comparison since it is particularly challenging for the reinitialization step, due to the high curvature of the interface in the tail region.

We observe that, at the coarsest level ($\ell_{\max}=8$), reinitialization increases both the mass and geometric errors compared to the non-reinitialized case. However, as the mesh is refined, the reinitialized results rapidly become more accurate, achieving much smaller errors at the finest resolution ($\ell_{\max}=12$). In fact, reinitialization has two competing effects: it introduces an interface shift, while it also mitigates the errors that would otherwise accumulate during advection as the field deforms.
As explained in Section~\ref{sec:reinit}, the interface-shift error scales with the accuracy of the finite-element representation of the level-set, whereas the advection-error mitigation depends on how the field deforms during transport (and thus on several factors).
In our tests, the first contribution decreases faster with mesh refinement, making reinitialization beneficial on finer meshes. The resolution threshold depends on the flow configuration.

Notice that the convergence orders $p_{\mathrm{m}}$ and $p_{\mathrm{g}}$ are consistently higher when reinitialization is applied, for all the values of $\ell_{\max}$. In particular, when reinitialization is performed, $\bar p_{\mathrm{m}} \approx 3.41$ and $ \bar p_{\mathrm{g}} \approx  2.62$.


\begin{table}[htbp]
  \centering
  \caption{Vortex kinematic transport test (Quad9): mass and geometric errors ($E_{\mathrm{m}}$, $E_{\mathrm{g}}$) and convergence orders $p_{\mathrm{m}}$, $p_{\mathrm{g}}$ without and with level-set reinitialization.}
  \label{tab:vortex_no_vs_reinit}
  \begingroup
  \setlength{\tabcolsep}{5pt}
  \begin{tabular}{c cccc cccc}
    \toprule
    & \multicolumn{4}{c}{Without reinitialization}
    & \multicolumn{4}{c}{With reinitialization} \\
    \cmidrule(lr){2-5}\cmidrule(lr){6-9}
    $\ell_{\max}$
      & $E_{\mathrm{m}}$ & $p_{\mathrm{m}}$ & $E_{\mathrm{g}}$ & $p_{\mathrm{g}}$
      & $E_{\mathrm{m}}$ & $p_{\mathrm{m}}$ & $E_{\mathrm{g}}$ & $p_{\mathrm{g}}$ \\
    \midrule
    8
      & $4.64\mathrm{E}{-04}$ & --   
      & $6.39\mathrm{E}{-04}$ & --  
      & $8.33\mathrm{E}{-03}$ & --   
      & $1.02\mathrm{E}{-03}$ & -- \\
    9
      & $2.12\mathrm{E}{-04}$ & 1.13 
      & $8.67\mathrm{E}{-05}$ & 2.88
      & $8.18\mathrm{E}{-04}$ & 3.35 
      & $1.75\mathrm{E}{-04}$ & 2.54 \\
    10
      & $1.35\mathrm{E}{-04}$ & 0.65 
      & $1.82\mathrm{E}{-05}$ & 2.25
      & $1.15\mathrm{E}{-04}$ & 2.83 
      & $2.31\mathrm{E}{-05}$ & 2.92 \\
    11
      & $6.12\mathrm{E}{-05}$ & 1.14 
      & $5.66\mathrm{E}{-06}$ & 1.68
      & $8.92\mathrm{E}{-06}$ & 3.69 
      & $3.94\mathrm{E}{-06}$ & 2.55 \\
    12
      & $1.58\mathrm{E}{-05}$ & 1.96 
      & $1.97\mathrm{E}{-06}$ & 1.52
      & $6.36\mathrm{E}{-07}$ & 3.81 
      & $7.06\mathrm{E}{-07}$ & 2.48 \\
    \bottomrule
  \end{tabular}
  \endgroup
\end{table}

\begin{figure}[!htb]
\centering
\includegraphics[width=0.3\linewidth]{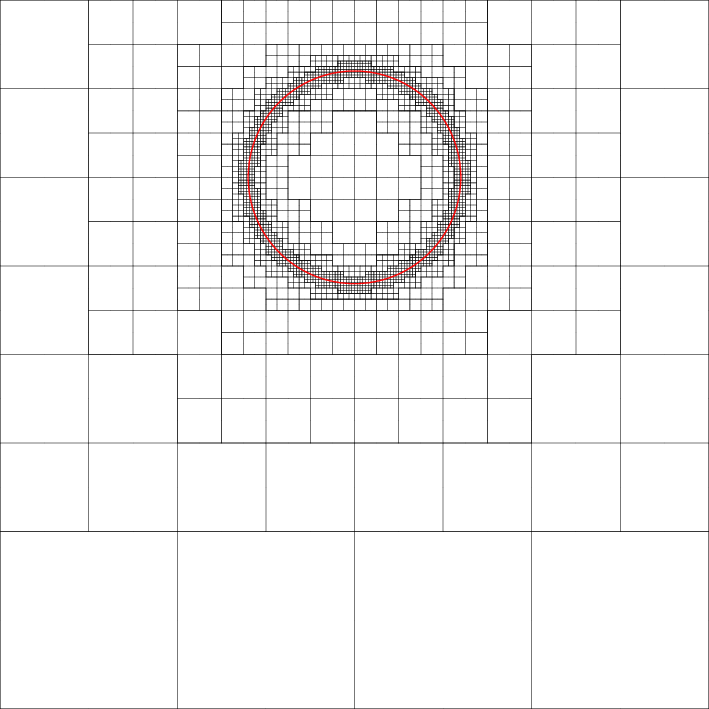}
\includegraphics[width=0.3\linewidth]{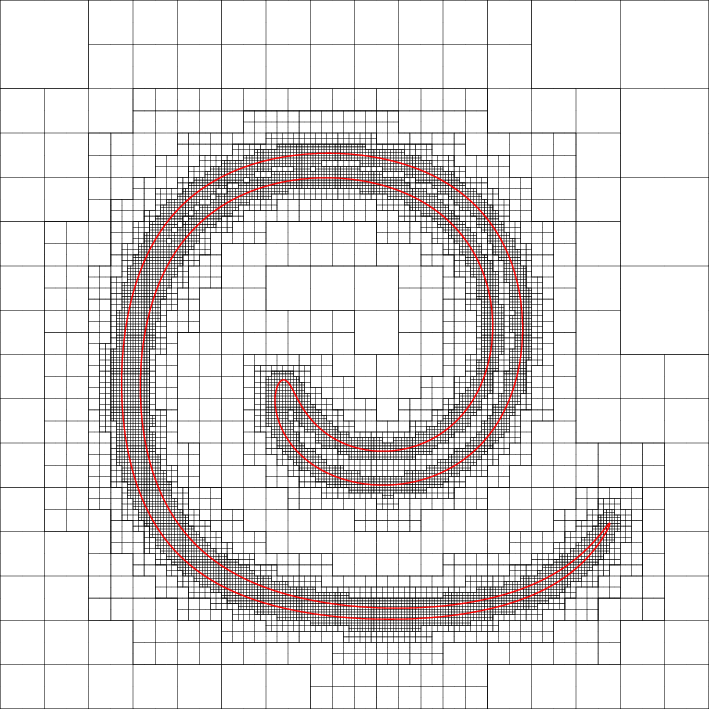}
\includegraphics[width=0.3\linewidth]{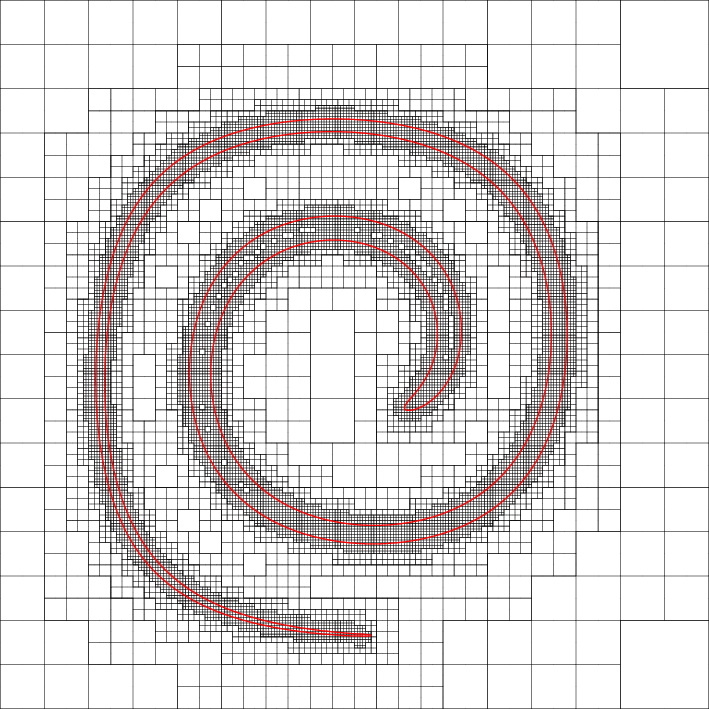}
\caption{Grid refinement for the 2D vortex test with Quad9 elements: $t=0$ on the left, $t=T/4$ in the center and $t=T/2$ on the right.}
\label{fig:init_2d}
\end{figure}

In Figure \ref{fig:init_2d}, the grid refinement for the 2D vortex test is reported, with the zero level-set curve depicted in red at different time instants, using a grid refinement of $\ell_{\max}=8$. Specifically, we report the initial grid configuration and the mesh refinement at half the period, when the interface undergoes its maximum stretching. The final time configuration is omitted, as the interface returns to its initial circular geometry.

In all subsequent simulations, the level-set function is reinitialized with
\(\tau_{\mathrm{reinit}} = 0.25\).

\paragraph{Simplicial elements: Tri7 and Tet15}

Table~\ref{tab:vortex_tri7_tet15} reports numerical results for the vortex benchmark obtained using simplicial finite elements, namely Tri7 in two dimensions and Tet15 in three dimensions. As in the Quad9 case, the level-set reinitialization is enabled for both discretizations. Results are reported as a function of the maximum refinement level $\ell_{\max}$, using the refinement ranges introduced above for the two- and three-dimensional cases: $\ell=8$--$12$ for Tri7 and $\ell=6$--$9$ for Tet15.

In two dimensions, the Tri7 discretization exhibits a clear reduction of both mass and geometric errors under mesh refinement. 
The observed average convergence rates are approximately $\bar{p}_{\mathrm{m}}\approx 3.55$ for the mass error and $\bar{p}_{\mathrm{g}}\approx 2.62$ for the geometric error, in line with the expected high-order behavior of the method.

In three dimensions, the Tet15 discretization shows a consistent decay of both error measures as the grid is refined. 
The corresponding average convergence rates are $\bar{p}_{\mathrm{m}}\approx 2.57$ and $\bar{p}_{\mathrm{g}}\approx 2.26$, which confirms robust accuracy also for simplicial elements in three dimensions.

\begin{table}[htbp]
  \centering
  \caption{Vortex kinematic transport test with level-set reinitialization:
           mass and geometric errors ($E_{\mathrm{m}}$, $E_{\mathrm{g}}$) and
           convergence orders ($p_{\mathrm{m}}$, $p_{\mathrm{g}}$) for
           Tri7 (2D) and Tet15 (3D) elements.}
  \label{tab:vortex_tri7_tet15}
  \begingroup
  \setlength{\tabcolsep}{5pt}
  \begin{tabular}{c cccc c cccc}
    \toprule
    & \multicolumn{4}{c}{Tri7} 
    & 
    & \multicolumn{4}{c}{Tet15} \\
    \cmidrule(lr){2-5}\cmidrule(lr){7-10}
    $\ell_{\max}$ 
      & $E_{\mathrm{m}}$ & $p_{\mathrm{m}}$ & $E_{\mathrm{g}}$ & $p_{\mathrm{g}}$
      & $\ell_{\max}$
      & $E_{\mathrm{m}}$ & $p_{\mathrm{m}}$ & $E_{\mathrm{g}}$ & $p_{\mathrm{g}}$ \\
    \midrule
    8  
      & $5.98\mathrm{E}{-03}$ & --    
      & $7.53\mathrm{E}{-04}$ & --    
      & 6  
      & $9.12\mathrm{E}{-02}$ & --    
      & $2.06\mathrm{E}{-03}$ & -- \\
    9  
      & $5.83\mathrm{E}{-04}$ & 3.36  
      & $1.26\mathrm{E}{-04}$ & 2.58  
      & 7  
      & $1.79\mathrm{E}{-02}$ & 2.35  
      & $4.80\mathrm{E}{-04}$ & 2.10 \\
    10 
      & $1.09\mathrm{E}{-04}$ & 2.42  
      & $1.78\mathrm{E}{-05}$ & 2.83  
      & 8  
      & $1.76\mathrm{E}{-03}$ & 3.34  
      & $8.75\mathrm{E}{-05}$ & 2.46 \\
    11 
      & $8.62\mathrm{E}{-06}$ & 3.66  
      & $2.82\mathrm{E}{-06}$ & 2.66  
      & 9  
      & $4.31\mathrm{E}{-04}$ & 2.03  
      & $1.87\mathrm{E}{-05}$ & 2.22 \\
    12 
      & $3.16\mathrm{E}{-07}$ & 4.77  
      & $5.24\mathrm{E}{-07}$ & 2.43  
      &     
      &                         &       
      &                         &      \\
    \bottomrule
  \end{tabular}
  \endgroup
\end{table}

\begin{figure}[!htb]
\centering
\includegraphics[width=0.45\linewidth]{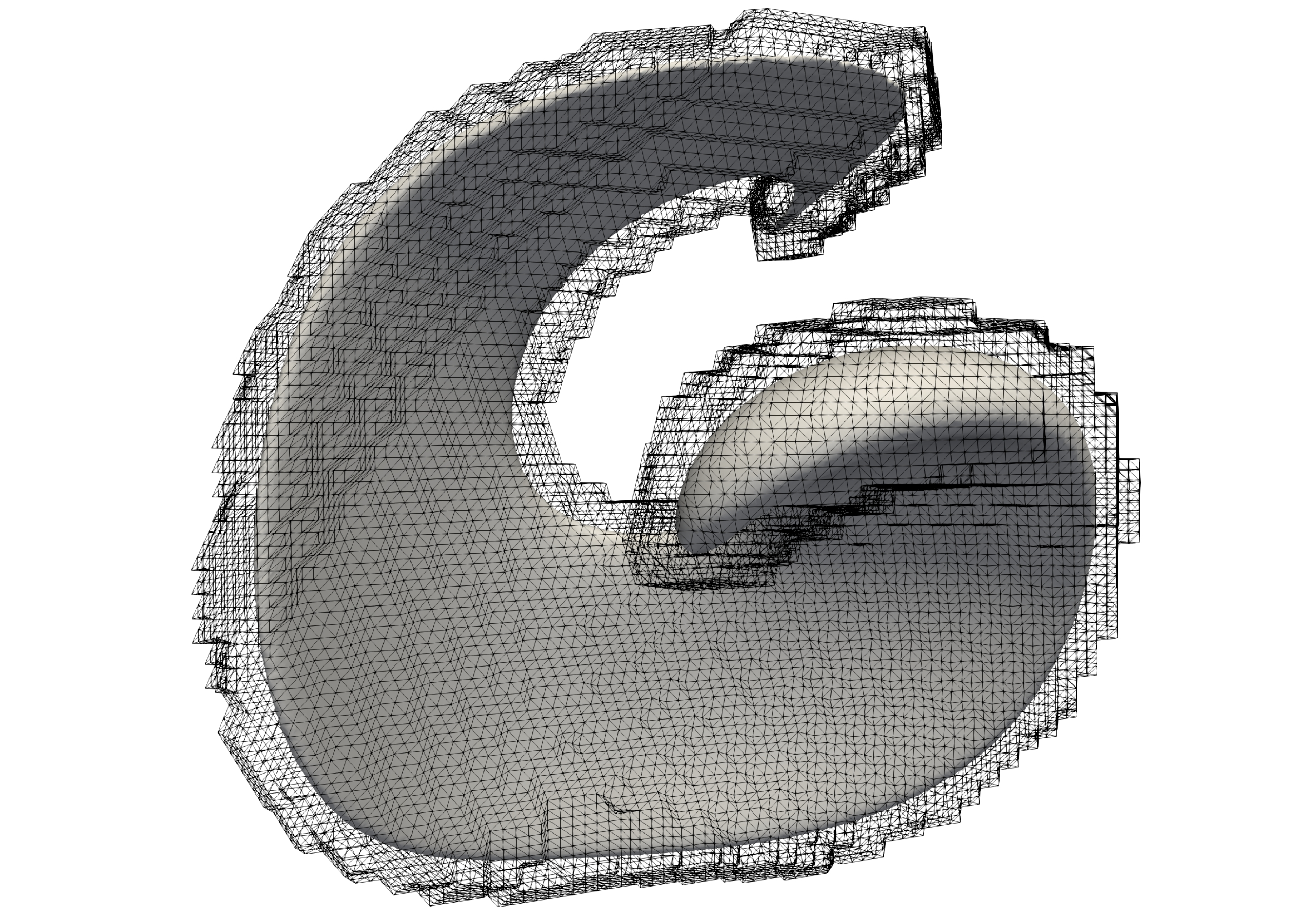}
\includegraphics[width=0.45\linewidth]{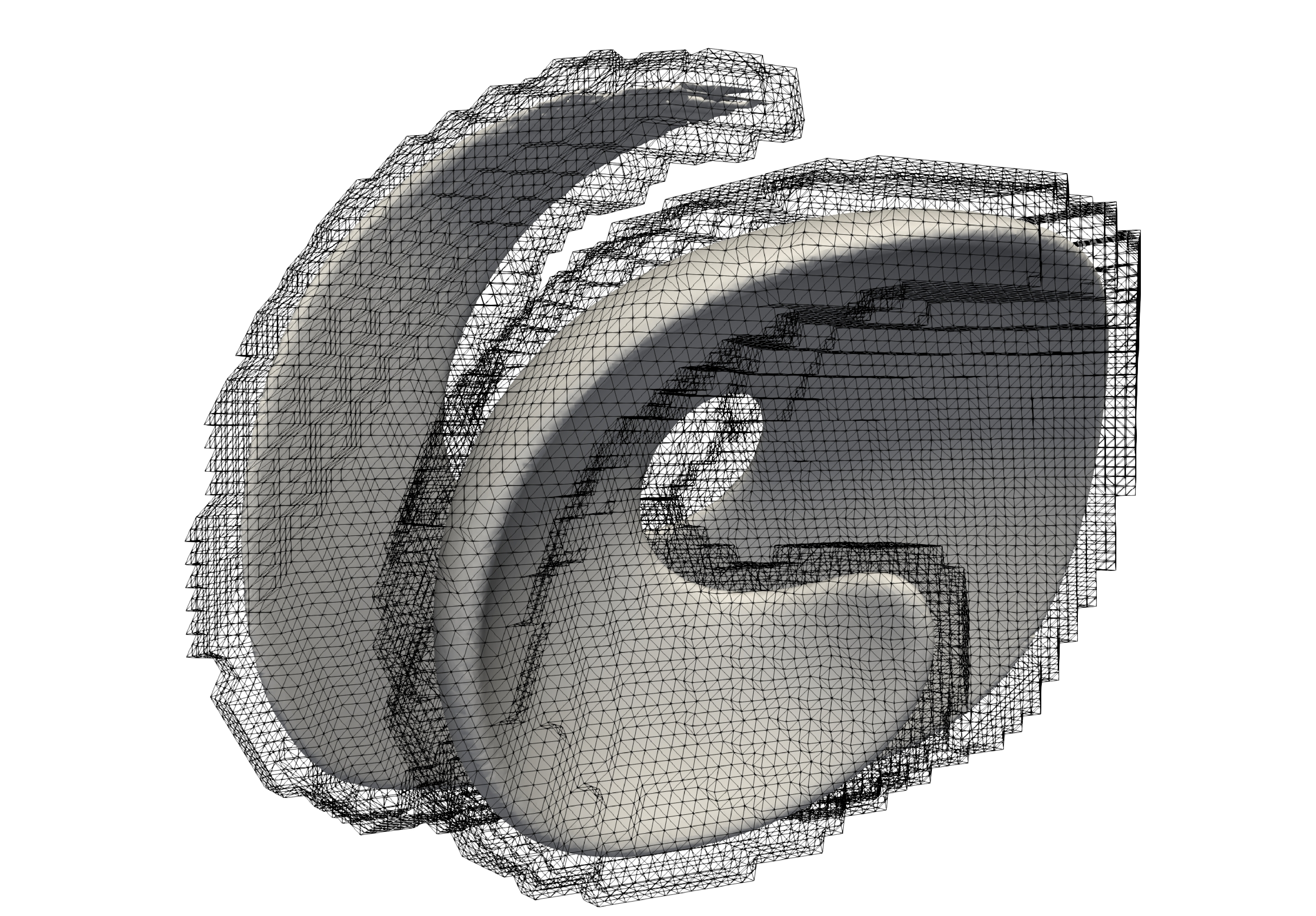}
\caption{Three-dimensional vortex test.
Interface deformation in the box domain computed with Tet15 elements at \(\ell_{\max}=7\).
Left: \(t=T/4\); right: \(t=T/2\).
Wireframe shows finest-level elements near the interface.}
\label{fig:v_3d}
\end{figure}

Figure~\ref{fig:v_3d} illustrates the interface evolution for the three-dimensional vortex test at two representative time instants.
The configuration at \(t = T/4\) corresponds to the onset of strong stretching, while the configuration at \(t = T/2\) exhibits extreme filamentation.
Results are shown for the Tet15 adaptive discretization with maximum refinement level \(\ell_{\max}=7\).
A wireframe representation of the mesh elements at the finest refinement level is also displayed in the vicinity of the interface.
Away from the interface, the adaptive mesh refinement procedure generates a graded mesh, which is not visible in the figure.

\paragraph{High-order 3D elements: Wedge21 and Hex27}

Table~\ref{tab:vortex_wedge21_hex27} reports numerical results for the vortex benchmark obtained using high-order three-dimensional finite elements, namely Wedge21 and Hex27. Both discretizations are tested over the same range of refinement levels, $\ell_{\max}=6$--$9$, enabling a direct comparison of accuracy and convergence behavior under identical resolution conditions. As in the previous cases, level-set reinitialization is enabled for both element types.



For the Wedge21 discretization, both the mass and geometric errors decrease consistently from $\ell=6$ to $\ell=9$, with average convergence rates $\bar{p}_{\mathrm{m}}\approx 2.46$ and $\bar{p}_{\mathrm{g}}\approx 2.33$, indicating robust high-order accuracy on prismatic elements.

For the Hex27 discretization, a similar decay is observed, with average convergence rates $\bar{p}_{\mathrm{m}}\approx 2.75$ and $\bar{p}_{\mathrm{g}}\approx 2.27$, confirming comparable accuracy on both hexahedral and wedge-based high-order elements.

Overall, the two discretizations exhibit very similar convergence trends for both mass conservation and geometric fidelity. Minor differences in the absolute error levels at coarse resolutions can be attributed to element topology and mapping distortions, while at finer levels, both Wedge21 and Hex27 achieve comparable accuracy, demonstrating the robustness of the method across different high-order 3D element families.

\begin{table}[htbp]
  \centering
  \caption{Vortex kinematic transport test with level-set reinitialization:
           mass and geometric errors ($E_{\mathrm{m}}$, $E_{\mathrm{g}}$) and
           convergence orders ($p_{\mathrm{m}}$, $p_{\mathrm{g}}$) for
           Wedge21 and Hex27 elements.}
  \label{tab:vortex_wedge21_hex27}
  \begingroup
  \setlength{\tabcolsep}{5pt}
  \begin{tabular}{c cccc cccc}
    \toprule
    & \multicolumn{4}{c}{Wedge21} 
    & \multicolumn{4}{c}{Hex27} \\
    \cmidrule(lr){2-5}\cmidrule(lr){6-9}
    $\ell_{\max}$
      & $E_{\mathrm{m}}$ & $p_{\mathrm{m}}$ & $E_{\mathrm{g}}$ & $p_{\mathrm{g}}$
      & $E_{\mathrm{m}}$ & $p_{\mathrm{m}}$ & $E_{\mathrm{g}}$ & $p_{\mathrm{g}}$ \\
    \midrule
    6  
      & $7.48\mathrm{E}{-02}$ & --    
      & $1.85\mathrm{E}{-03}$ & --    
      & $1.45\mathrm{E}{-01}$ & --    
      & $3.03\mathrm{E}{-03}$ & -- \\
    7  
      & $1.81\mathrm{E}{-02}$ & 2.05  
      & $4.75\mathrm{E}{-04}$ & 1.96  
      & $2.55\mathrm{E}{-02}$ & 2.50  
      & $6.51\mathrm{E}{-04}$ & 2.22 \\
    8  
      & $1.98\mathrm{E}{-03}$ & 3.20  
      & $7.96\mathrm{E}{-05}$ & 2.58  
      & $2.51\mathrm{E}{-03}$ & 3.35  
      & $1.34\mathrm{E}{-04}$ & 2.28 \\
    9  
      & $4.48\mathrm{E}{-04}$ & 2.14  
      & $1.47\mathrm{E}{-05}$ & 2.44  
      & $4.75\mathrm{E}{-04}$ & 2.40  
      & $2.68\mathrm{E}{-05}$ & 2.32 \\
    \bottomrule
  \end{tabular}
  \endgroup
\end{table}


\paragraph{Comparison with literature}

To further assess the proposed AMR strategy, we compare the mass conservation properties with the results reported by Ngo and Choi~\cite{ngo2017multi}, who considered the same two-dimensional time-reversed single-vortex test with period $T=8$ using an adaptive finite-element level-set method. From the number of elements of their uniformly refined reference mesh, the minimum local resolution near the interface can be estimated as $h_{\min}\approx 1/512$ using triangular elements. In our framework, this resolution corresponds to refinement level $\ell=9$, considering biquadratic triangular elements (Tri7).

At this resolution, the relative mass error obtained in the present work is 
$E_\text{m} = 5.83\times 10^{-4}$ (i.e.\ $0.058\%$), which is significantly lower than the mass errors reported in~\cite{ngo2017multi} for both one-level and two-level adaptive meshes, as shown in Table~\ref{tab:vortex_literature}. Despite differences in discretization and AMR strategies, the comparison is meaningful at the level of local interface resolution. Ngo and Choi~\cite{ngo2017multi} use linear triangular elements, resulting in fewer global degrees of freedom, while their one-level (M1) and two-level (M2) adaptive meshes preserve a similar minimum mesh spacing near the interface and differ mainly in the extent of the refined region. Accordingly, the present results are compared at matched $h_{\min}$ rather than matched degrees of freedom, confirming competitive and improved mass conservation.
\begin{table}[h!]
\centering
\caption{Comparison of relative mass error for the 2D vortex test at comparable interface resolution $h_{\min}\approx 1/512$. 
Numbers of nodes and elements are evaluated at $t=T/2$ (maximum stretching).}
\label{tab:vortex_literature}
\begin{tabular}{lcccc}
\toprule
Method & \# nodes & \# elements & $h_{\min}$ & Mass error (\%) \\
\midrule
This work (AMR, Tri7, $\ell=9$) & 168\,491 & 54\,467 & $1/512$ & 0.058 \\
Ngo \& Choi (2016), M1 & 82\,391 & 163\,824 & $\approx 1/512$ & 0.231 \\
Ngo \& Choi (2016), M2 & 44\,056 & 87\,634 & $\approx 1/512$ & 0.448 \\
\bottomrule
\end{tabular}
\end{table}

\subsection{Rising-bubble-type test}
\label{sec:rb}

\paragraph{Test definition and funnel geometry}

The rising-bubble-type kinematic test is driven by a prescribed divergence-free, time-periodic velocity field on two- and three-dimensional domains.
In two dimensions, we define
\begin{equation}
{v}(x,y,t)
=
A\left\langle -U(x)\,V'(y),\; U'(x)\,V(y)\right\rangle
\cos\!\left(\frac{\pi t}{T}\right),
\end{equation}
where \(A=0.1\) and \(T=8\).
The auxiliary functions are
\begin{equation}
U(s)=\cos^2\!\left(\frac{\pi s}{L}\right)\sin\!\left(\frac{2\pi s}{L}\right),
\qquad
V(s)=\sin^2\!\left(\frac{\pi s}{L}\right),
\qquad L=1,
\end{equation}
with derivatives
\begin{equation}
U'(s)=\frac{\pi}{L}\left[\cos\!\left(\frac{2\pi s}{L}\right)+\cos\!\left(\frac{4\pi s}{L}\right)\right],
\qquad
V'(s)=\frac{\pi}{L}\sin\!\left(\frac{2\pi s}{L}\right).
\end{equation}

In three dimensions, we set \(T=4\) and define
\begin{equation}
{v}(x,y,z,t)
=
A\left\langle
-\,U(x)\,V'(z),\;
-\,U(y)\,V'(z),\;
\big(U'(x)+U'(y)\big)\,V(z)
\right\rangle
\cos\!\left(\frac{\pi t}{T}\right),
\end{equation}
with the same functions \(U,V\) (and \(L=1\)).

The initial interface is a circular bubble of radius $r=0.15$ centered at $(0,0.25)$ in two dimensions, and a spherical bubble of the same radius centered at $(0,0,0.25)$ in three dimensions. During half a period, the interface rises and elongates upward in a motion reminiscent of buoyant ascent, then returns to its original configuration at $t=T$. Unlike the vortex benchmark in Section~\ref{sec:vortex}, where the level-set is stretched along the flow direction, here it is compressed. In fact, the velocity-field gradient is opposite to the flow direction. Both flows deform the level-set in the flow direction, but in two opposite ways. Because advection responds differently to stretching than to compression, this test is designed to evaluate how the proposed approach performs in this setting. A more thorough discussion of how level-set deformation influences advection accuracy is provided in Section~\ref{sec:reinit}.

This benchmark is run directly on a \emph{funnel-shaped} grid in both two and three dimensions, as shown in Figure~\ref{fig:grids}. In two dimensions, the lateral boundaries are defined as hyperbolic arcs:
\[
x(y) = h \pm a\sqrt{1 + \left(\frac{y - k}{b}\right)^2},
\]
with $a=0.4$ and $b=0.5$ controlling the aperture and curvature, and $(h,k) = (0.0, 0.5)$ denoting the funnel axis position. 
In three dimensions, the funnel geometry is obtained by smoothly remapping the Cartesian cube mesh into a funnel-shaped domain.
The transformation scales the lateral coordinates as a function of the vertical position, so that each horizontal cross-section remains square while its side length varies along the height according to the prescribed funnel profile.
The resulting meshes therefore exhibit a gradual lateral contraction or expansion, providing a controlled amount of element distortion.

\paragraph{Element-type comparison: Tri7 and Wedge21 on funnel meshes}
Table~\ref{tab:rising_bubble_funnel} reports the mass and geometric errors, together with the corresponding convergence orders, for the rising-bubble test computed on funnel-shaped grids using Tri7 elements in two dimensions and Wedge21 elements in three dimensions. In both cases, level-set reinitialization is enabled. Results are reported as a function of the maximum refinement level $\ell_{\max}$ (or equivalently of the minimum local spacing $h_{\min}$).
\begin{table}[htbp]
  \centering
  \caption{Rising-bubble-type kinematic transport test on funnel meshes with level-set reinitialization:
           mass and geometric errors ($E_{\mathrm{m}}$, $E_{\mathrm{g}}$) and
           convergence orders ($p_{\mathrm{m}}$, $p_{\mathrm{g}}$) for Tri7 (2D) and Wedge21 (3D) elements.}
  \label{tab:rising_bubble_funnel}
  \begingroup
  \setlength{\tabcolsep}{5pt}
  \begin{tabular}{c cccc c cccc}
    \toprule
    & \multicolumn{4}{c}{Tri7} 
    & 
    & \multicolumn{4}{c}{Wedge21} \\
    \cmidrule(lr){2-5}\cmidrule(lr){7-10}
    $\ell_{\max}$
      & $E_{\mathrm{m}}$ & $p_{\mathrm{m}}$ & $E_{\mathrm{g}}$ & $p_{\mathrm{g}}$
      & $\ell_{\max}$
      & $E_{\mathrm{m}}$ & $p_{\mathrm{m}}$ & $E_{\mathrm{g}}$ & $p_{\mathrm{g}}$ \\
    \midrule
    8  & $1.51\mathrm{E}{-01}$ & --    
       & $1.07\mathrm{E}{-02}$ & --    
       & 6  & $1.18\mathrm{E}{-02}$ & --    
            & $3.27\mathrm{E}{-04}$ & -- \\
    9  & $7.92\mathrm{E}{-05}$ & 10.89 
       & $9.43\mathrm{E}{-06}$ & 10.14 
       & 7  & $1.94\mathrm{E}{-03}$ & 2.60  
            & $6.33\mathrm{E}{-05}$ & 2.37 \\
    10 & $6.06\mathrm{E}{-06}$ & 3.71  
       & $9.90\mathrm{E}{-07}$ & 3.25  
       & 8  & $1.36\mathrm{E}{-04}$ & 3.84  
            & $9.20\mathrm{E}{-06}$ & 2.78 \\
    11 & $3.63\mathrm{E}{-07}$ & 4.06  
       & $8.95\mathrm{E}{-08}$ & 3.47  
       & 9  & $5.53\mathrm{E}{-06}$ & 4.62  
            & $1.26\mathrm{E}{-06}$ & 2.87 \\
    12 & $4.10\mathrm{E}{-08}$ & 3.15  
       & $1.07\mathrm{E}{-08}$ & 3.07  
       &     &                         &       
            &                         &      \\
    \bottomrule
  \end{tabular}
  \endgroup
\end{table}

The results reported in Table~\ref{tab:rising_bubble_funnel} show that both Tri7 and Wedge21 elements achieve monotonic convergence of the mass and geometric errors under compressive deformation on funnel-shaped meshes. In two dimensions, the Tri7 discretization reaches mass and geometric errors of the order of $10^{-7}$ and $10^{-9}$, respectively, at the finest refinement level, while in three dimensions the Wedge21 discretization attains $E_{\mathrm{m}}\approx 5\times 10^{-6}$ and $E_{\mathrm{g}}\approx 10^{-6}$ at $\ell_{\max}=9$. Despite the strongly anisotropic deformation and the curved element topology induced by the funnel geometry, the observed convergence rates remain close to the theoretical expectations, with approximately third-order accuracy for the mass error and second-order accuracy for the geometric error. These results confirm the robustness of the proposed AMR level-set formulation under compression-dominated interface transport and distorted mesh configurations. 

The level-8 Tri7 errors are noticeably larger than those observed for the other configurations.
This behavior is attributed to the fact that, at this resolution, the deformed bubble approaches the funnel boundaries.
The imposed level-set boundary condition, \(\PsiLS=-1\), on walls crossed by the velocity field is unphysical in this regime and leads to the formation of spurious bubbles.
A more appropriate outflow boundary condition should therefore be employed.
Nevertheless, the algorithm remains robust and does not exhibit any numerical breakdown.

\begin{figure}[!htb]
\centering
\includegraphics[width=0.3\linewidth]{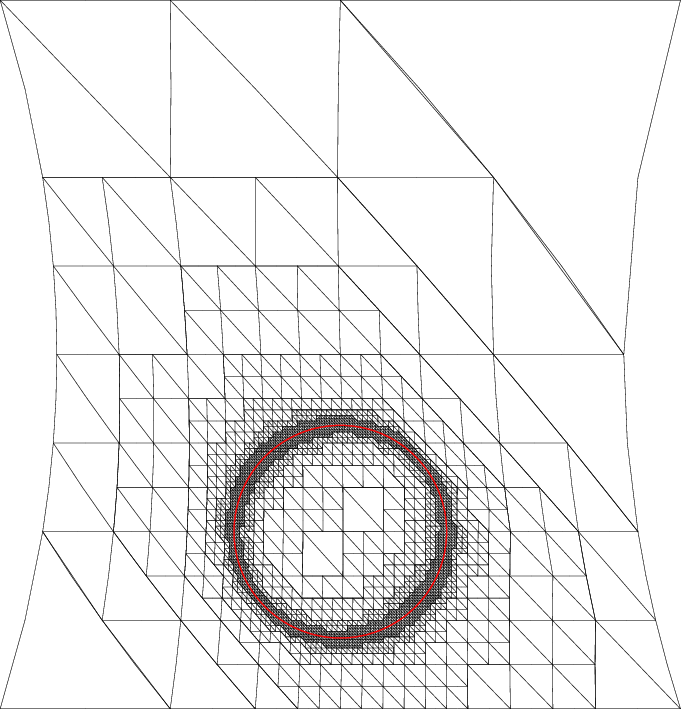}
\includegraphics[width=0.3\linewidth]{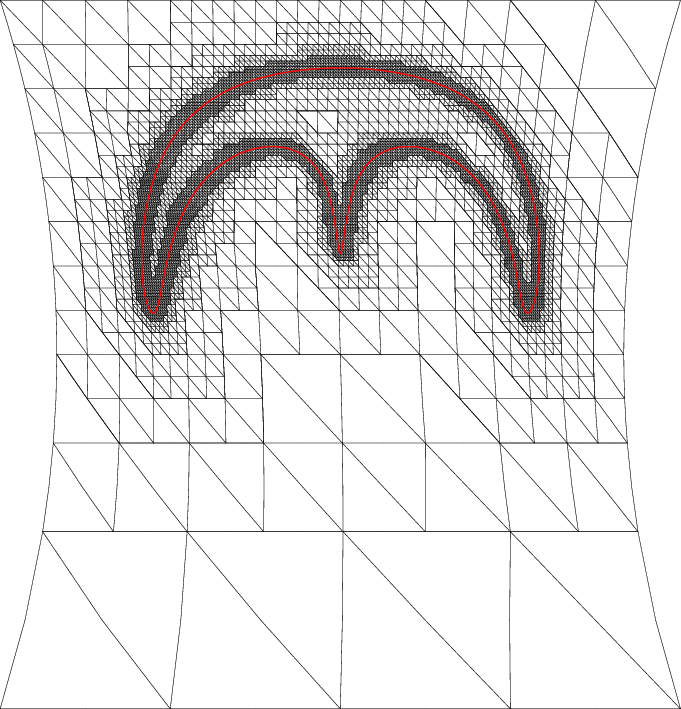}
\includegraphics[width=0.3\linewidth]{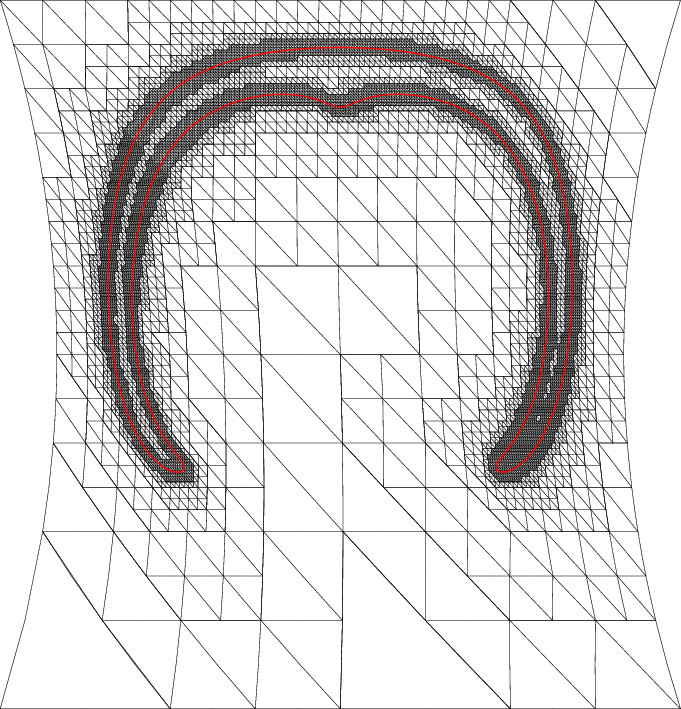}
\caption{Grid refinement for the 2D rising bubble test with Tri7 elements at $\ell_{\max} = 8$: $t=0$ on the left, $t=T/4$ in the center and $t=T/2$ on the right.}
\label{fig:init_2d_funnel}
\end{figure}

In Figure~\ref{fig:init_2d_funnel}, we report, for the two-dimensional funnel case, the interface position of the bubble at four different time steps of the simulation, considering a grid refinement with $\ell_{\max}=8$. We notice that the initial interface is rising to the top of the domain, where it starts to stretch along the $x$-direction, until the tails are moved towards the bottom of the domain. After that, the interface returns back to the initial configuration, where the circular geometry is achieved with low error.

\begin{figure}[!htb]
\centering
\includegraphics[width=0.45\linewidth]{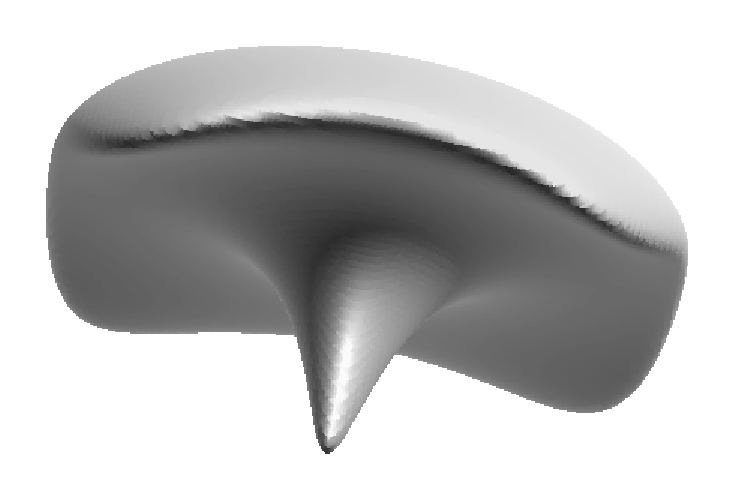}
\includegraphics[width=0.45\linewidth]{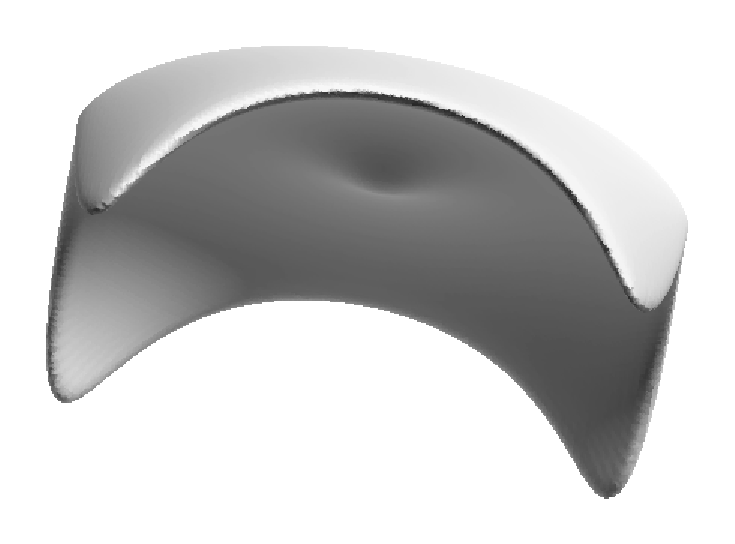}
\caption{Rising bubble 3D test interface deformation at different time steps for the funnel domain with Wedge21 elements at $\ell_{\max}=9$: left $t=T/4$, and right $t=T/2$.}
\label{fig:rb_3d}
\end{figure}

In Figure~\ref{fig:rb_3d}, the three-dimensional interface evolution is shown at two representative times, considering the half-period configuration with maximum elongation inside the funnel. 

\subsection{Rigid-body rotation test}
\label{sec:rot}

\paragraph{Test definition}
The rotational test, originally introduced by LeVeque~\cite{leveque1996high} and widely adopted in subsequent studies~\cite{Rider1998,Enright2002,Aulisa2003Mixed,Henri2022Geometrical}, serves as a baseline for verifying the consistency of the advection scheme in a purely rigid-body motion.

In two dimensions, the velocity field is defined as
\begin{equation}
{v}(x,y) = (y,\,-x),
\end{equation}
which corresponds to a clockwise rotation about the origin. The computational domain is the unit square $\Omega=[-0.5,0.5]^2$, discretized using Tri7 elements. The initial level-set field $\PsiLS(x,y)$ is initialized as the signed distance from a circle of radius $r=0.15$ centered at $(0.0,0.25)$. During the simulation, the interface undergoes a rigid rotation of $2\pi$ radians, theoretically returning to its initial position after one full period.

In three dimensions, we consider a fully three-dimensional rotational field inducing rotation around the main diagonal of the computational domain, defined as
\begin{equation}
{v}(x,y,z)=(y-z,\; z-x,\; x-y)\,.
\end{equation}
The computational domain is the \emph{ball-shaped} centered at the origin with radius equal to 0.5, discretized using Tet15 elements, which introduces curved simplicial elements. The initial interface is a sphere of radius $r=0.15$ centered at $(0.0,0.0,0.25)$, initialized as a signed-distance function. As in the two-dimensional case, the interface undergoes a rigid-body motion and should ideally recover its initial configuration after one full period.

For ideal rigid rotation, both the mass error $E_{\mathrm{m}}$ and the geometric error $E_{\mathrm{g}}$ remain negligible and converge with the formal accuracy of the advection method. Since the prescribed velocity fields induce a rigid-body motion of the level-set, the interface is transported without deformation and does not require reinitialization; this behavior is automatically enforced by the reinitialization criterion described in Section~\ref{sec:reinit} and is shown quantitatively in Section~\ref{sec:reinit_perf}.

\paragraph{Error convergence on Tri7 (square) and Tet15 (ball) domains}
Table~\ref{tab:rotation_tri7_tet15} reports the mass and geometric errors, together with the corresponding convergence orders, for the rigid-body rotation test computed on a straight square domain using Tri7 elements in two dimensions and on a curved spherical domain using Tet15 elements in three dimensions. Results are reported as a function of the maximum refinement level $\ell_{\max}$ (or equivalently of the minimum local spacing $h_{\min}$).

In two dimensions, the Tri7 discretization yields extremely small errors, reaching $E_{\mathrm{m}}\approx 8\times 10^{-9}$ and $E_{\mathrm{g}}\approx 1\times 10^{-9}$ at $\ell_{\max}=12$. The observed convergence rates are close to second order for both mass and geometric errors, in agreement with theoretical expectations for rigid-body transport. 

In three dimensions, the Tet15 discretization also produces small geometric errors, with $E_{\mathrm{g}}\approx 6\times 10^{-7}$ at $\ell_{\max}=9$. The mass error decreases initially with refinement but exhibits a mild saturation at fine levels, stabilizing around $E_{\mathrm{m}}\approx 3\times 10^{-5}$. This behavior is attributed to the combined effects of curved element geometry, high-order quadrature, and accumulated round-off errors. 

These results confirm that the proposed AMR level-set method preserves good accuracy under rigid-body motion on both straight and curved meshes, while correctly suppressing unnecessary reinitialization.

\begin{table}[htbp]
  \centering
  \caption{Rigid-body rotation kinematic transport test:
           mass and geometric errors ($E_{\mathrm{m}}$, $E_{\mathrm{g}}$) and
           convergence orders ($p_{\mathrm{m}}$, $p_{\mathrm{g}}$) for
           Tri7 (2D square) and Tet15 (3D ball) elements.}
  \label{tab:rotation_tri7_tet15}
  \begingroup
  \setlength{\tabcolsep}{5pt}
  \begin{tabular}{c cccc c cccc}
    \toprule
    & \multicolumn{4}{c}{Tri7} 
    & 
    & \multicolumn{4}{c}{Tet15} \\
    \cmidrule(lr){2-5}\cmidrule(lr){7-10}
    $\ell_{\max}$
      & $E_{\mathrm{m}}$ & $p_{\mathrm{m}}$ & $E_{\mathrm{g}}$ & $p_{\mathrm{g}}$
      & $\ell_{\max}$
      & $E_{\mathrm{m}}$ & $p_{\mathrm{m}}$ & $E_{\mathrm{g}}$ & $p_{\mathrm{g}}$ \\
    \midrule
    8  & $4.49\mathrm{E}{-06}$ & --    & $1.39\mathrm{E}{-06}$ & --    
       & 5  & $1.95\mathrm{E}{-04}$ & --    & $1.51\mathrm{E}{-04}$ & -- \\
    9  & $3.29\mathrm{E}{-07}$ & 3.77  & $8.80\mathrm{E}{-08}$ & 3.98  
       & 6  & $1.20\mathrm{E}{-04}$ & 0.69  & $9.96\mathrm{E}{-06}$ & 3.92 \\
    10 & $3.52\mathrm{E}{-07}$ & -0.10 & $2.49\mathrm{E}{-08}$ & 1.82  
       & 7  & $3.47\mathrm{E}{-05}$ & 1.79  & $1.21\mathrm{E}{-06}$ & 3.05 \\
    11 & $1.10\mathrm{E}{-07}$ & 1.67  & $7.81\mathrm{E}{-09}$ & 1.67  
       & 8  & $2.95\mathrm{E}{-05}$ & 0.23  & $6.84\mathrm{E}{-07}$ & 0.82 \\
    12 & $8.01\mathrm{E}{-09}$ & 3.79  & $1.33\mathrm{E}{-09}$ & 2.56  
       & 9  & $2.93\mathrm{E}{-05}$ & 0.01  & $5.93\mathrm{E}{-07}$ & 0.20 \\
    \bottomrule
  \end{tabular}
  \endgroup
\end{table}



\subsection{Reinitialization frequency}
\label{sec:reinit_perf}

In this section, we analyze both the frequency and the computational cost of the reinitialization procedure.
Unlike approaches based on a fixed reinitialization interval, the criterion employed in this work (see Sec.~\ref{sec:reinit}) activates reinitialization dynamically based on the evolution of the level-set field and on the prescribed threshold parameter~$\overline{D}$.
As a result, the total number of reinitializations is not imposed \emph{a priori}, but emerges naturally from the simulation. On the other hand, this approach requires evaluating the global indicator $\overline{D}$ at each time step. If the indicator exceeds the threshold $\tau_{\mathrm{reinit}}$, the reinitialization procedure is triggered. This means that the computational cost includes a small cost at every time step to compute the indicator, plus an additional cost only when the reinitialization is activated.

Table~\ref{tab:number_reinit} reports the total number of reinitializations observed in the test cases considered in this work, for the different element types and benchmarks.  
For compactness, the following acronyms are used:  
\emph{VQ9} = 2D vortex test with Quad9 elements;  
\emph{VT7} = 2D vortex test with Tri7 elements;  
\emph{RB-T7} = 2D rising-bubble test on funnel meshes with Tri7 elements;  
\emph{V-T15} = 3D vortex test with Tet15 elements;  
\emph{V-W21} = 3D vortex test with Wedge21 elements;  
\emph{V-H27} = 3D vortex test with Hex27 elements;  
\emph{RB-W21} = 3D rising-bubble test on funnel meshes with Wedge21 elements.

For the vortex benchmark, the number of reinitializations is essentially grid-independent in two dimensions (VQ9 and VT7) and becomes nearly grid-independent in three dimensions (V-T15, V-W21 and V-H27), where only minor variations of one or two events are observed across refinement levels.

For the rising-bubble-type benchmark computed directly on funnel meshes, the number of level-set reinitializations stabilizes rapidly in three dimensions (RB--W21).
In two dimensions (RB--T7), an elevated number of reinitializations is observed at the level-8 refinement level.
This effect is consistent with the boundary-related artifacts discussed above.
Despite this behavior, the method remains stable under the same non-optimal boundary treatment.
Upon mesh refinement, the number of reinitializations decreases and aligns with the values observed in the vortex test cases.

In the rotation test cases, both in two and three dimensions, no reinitialization is triggered. Since the prescribed velocity field induces a rigid-body motion of the level-set, the only deformation arises from numerical errors in the advection step. For the chosen value of $\tau_{\mathrm{reinit}}$, these errors are not sufficient to activate the reinitialization procedure. 

\begin{table}[htbp]
\centering
\caption{Total number of reinitializations for the benchmarks considered in this work. Acronyms are defined in the text.}
\label{tab:number_reinit}
\renewcommand{\arraystretch}{1.05}
\setlength{\tabcolsep}{4pt}
\begin{tabular}{c ccc c cccc}
\toprule
& \multicolumn{3}{c}{\textbf{2D}}
& 
& \multicolumn{4}{c}{\textbf{3D}} \\
\cmidrule(lr){2-4}\cmidrule(lr){6-9}
$\ell_{\max}$ 
& VQ9 & VT7 & RB-T7
&
$\ell_{\max}$
& V-T15 & V-W21 & V-H27 & RB-W21 \\
\midrule
8  & 15 & 15 & 89  & 6 & 12 & 10 & 11 & 10 \\
9  & 15 & 15 & 15  & 7 & 11 & 11 & 11 & 10 \\
10 & 15 & 15 & 15  & 8 & 11 & 11 & 11 & 10 \\
11 & 15 & 15 & 15  & 9 & 11 & 11 & 11 & 10 \\
12 & 15 & 15 & 15  &   &    &    &    &    \\
\bottomrule
\end{tabular}
\end{table}

\section{Conclusions} 

We have presented an adaptive finite element framework for level-set transport in two-phase flow simulations that combines the robustness of implicit interface representation with the accuracy of marker-based transport. The level-set formulation naturally accommodates topological changes, while forward marker transport enables predictive adaptive mesh construction and backward characteristic tracing provides accurate evaluation of the advected level-set field. An efficient multilevel point-location strategy supports both operations. A consistent reinitialization procedure restores the signed-distance property while preserving interface geometry and mass conservation.

Numerical experiments in two and three spatial dimensions, on structured and unstructured meshes and across multiple element types, demonstrate that the proposed adaptive strategy attains accuracy comparable to uniform discretizations at substantially reduced computational cost. The method exhibits robust behavior under strong deformation and aggressive adaptivity, while maintaining stable and conservative interface evolution.

The present study considers prescribed velocity fields to isolate transport and adaptivity mechanisms. The next step is coupling with multiphase flow solvers, enabling fully adaptive finite-element simulations of complex two-phase flows that combine topology-resolving level-set evolution with marker-enhanced accuracy. This establishes a marker-consistent, fully adaptive foundation for accurate and efficient level-set transport.

\bibliographystyle{elsarticle-num}
\bibliography{bib}

@article{OsherSethian1988,
  title   = {Fronts Propagating with Curvature-Dependent Speed: Algorithms Based on {Hamilton--Jacobi} Formulations},
  author  = {Osher, Stanley and Sethian, James A.},
  journal = {Journal of Computational Physics},
  volume  = {79},
  number  = {1},
  pages   = {12--49},
  year    = {1988}
}

@article{JiangPeng2000,
  title={Weighted {ENO} Schemes for {Hamilton--Jacobi} Equations},
  author={Jiang, Guang-Shan and Peng, Dongsu},
  journal={SIAM Journal on Scientific Computing},
  volume={21},
  number={6},
  pages={2126--2143},
  year={2000}
}

@article{AdalsteinssonSethian1995,
  title={A Fast Level Set Method for Propagating Interfaces},
  author={Adalsteinsson, David and Sethian, James A.},
  journal={Journal of Computational Physics},
  volume={118},
  number={2},
  pages={269--277},
  year={1995}
}

@article{MinGibou2007,
  title={A second order accurate level set method on non-graded adaptive Cartesian grids},
  author={Min, Chohong and Gibou, Fr{\'e}d{\'e}ric},
  journal={Journal of Computational Physics},
  volume={225},
  number={1},
  pages={300--321},
  year={2007},
  publisher={Elsevier}
}

@article{Popinet2003,
  title={Gerris: A Tree-Based Adaptive Solver for the Incompressible Euler Equations in Complex Geometries},
  author={Popinet, St{\'e}phane},
  journal={Journal of Computational Physics},
  volume={190},
  number={2},
  pages={572--600},
  year={2003}
}

@article{Shu2009,
  title={High Order Weighted Essentially Nonoscillatory Schemes for Convection Dominated Problems},
  author={Shu, Chi-Wang},
  journal={SIAM Review},
  volume={51},
  number={1},
  pages={82--126},
  year={2009}
}

@article{StaniforthCote1991,
  title={Semi-Lagrangian Integration Schemes for Atmospheric Models: A Review},
  author={Staniforth, Andrew and C{\^o}t{\'e}, Jean},
  journal={Monthly Weather Review},
  volume={119},
  number={9},
  pages={2206--2223},
  year={1991}
}

@article{Olsson2007,
  title={A Conservative Level Set Method for Two Phase Flow},
  author={Olsson, E. and Kreiss, G. and Zahedi, S.},
  journal={Journal of Computational Physics},
  volume={225},
  number={1},
  pages={785--807},
  year={2007}
}

@article{rider1998,
  title     = {Reconstructing volume tracking},
  author    = {W. J. Rider and D. B. Kothe},
  journal   = {Journal of Computational Physics},
  volume    = {141},
  number    = {2},
  pages     = {112--152},
  year      = {1998}
}

@article{enright2002,
  title     = {A hybrid particle level set method for improved interface capturing},
  author    = {D. Enright and R. Fedkiw and J. Ferziger and I. Mitchell},
  journal   = {Journal of Computational Physics},
  volume    = {183},
  number    = {1},
  pages     = {83--116},
  year      = {2002}
}

@article{aulisa2004surface,
  title={A surface marker algorithm coupled to an area-preserving marker redistribution method for three-dimensional interface tracking},
  author={Aulisa, Eugenio and Manservisi, Sandro and Scardovelli, Ruben},
  journal={Journal of Computational Physics},
  volume={197},
  number={2},
  pages={555--584},
  year={2004},
  publisher={Elsevier}
}

@article{aulisa2007interface,
  title={Interface reconstruction with least-squares fit and split advection in three-dimensional Cartesian geometry},
  author={Aulisa, Eugenio and Manservisi, Sandro and Scardovelli, Ruben and Zaleski, Stephane},
  journal={Journal of Computational Physics},
  volume={225},
  number={2},
  pages={2301--2319},
  year={2007},
  publisher={Elsevier}
}

@article{aulisa2025mixed,
  title={A Mixed Marker and Level-Set Front-Tracking Approximation for Multiphase Flows Simulations},
  author={Aulisa, Eugenio and Barbi, Giacomo and Chierici, Andrea and Manservisi, Sandro},
  journal={Journal of Computational Physics},
  pages={114285},
  year={2025},
  publisher={Elsevier}
}

@article{henri2022geometrical,
  title={Geometrical level set reinitialization using closest point method and kink detection for thin filaments, topology changes and two-phase flows},
  author={Henri, F{\'e}lix and Coquerelle, Mathieu and Lubin, Pierre},
  journal={Journal of Computational Physics},
  volume={448},
  pages={110704},
  year={2022},
  publisher={Elsevier}
}

@article{sussman1994level,
  title={A level set approach for computing solutions to incompressible two-phase flow},
  author={Sussman, Mark and Smereka, Peter and Osher, Stanley},
  journal={Journal of Computational physics},
  volume={114},
  number={1},
  pages={146--159},
  year={1994},
  publisher={Elsevier}
}

@article{aulisa2003mixed,
  title={A mixed markers and volume-of-fluid method for the reconstruction and advection of interfaces in two-phase and free-boundary flows},
  author={Aulisa, Eugenio and Manservisi, Sandro and Scardovelli, Ruben},
  journal={Journal of Computational Physics},
  volume={188},
  number={2},
  pages={611--639},
  year={2003},
  publisher={Elsevier}
}

@article{leveque1996high,
  title={High-resolution conservative algorithms for advection in incompressible flow},
  author={Leveque, Randall J},
  journal={SIAM Journal on Numerical Analysis},
  volume={33},
  number={2},
  pages={627--665},
  year={1996},
  publisher={SIAM}
}

@article{shakoor2015kdtree,
  title={An efficient and parallel level set reinitialization method – Application to micromechanics and microstructural evolutions},
  author={Shaakor, Modesar and Scholtes, Benjamin and Bouchard, Pierre-Olivier and Bernacki, Marc},
  journal={Applied Mathematical Modelling},
  volume={39},
  number={23-24},
  pages={7291-7302},
  year={2015},
  publisher={Elsevier}
}

@techreport{parolini2007reinit,
  title={A local projection reinitialization procedure for the level set equation on unstructured grids},
  author={Parolini, Nicola and Burman, Erik},
  year={2007},
  institution={Technical Report CMCS-REPORT-2007-004, {\'E}cole Polytechnique F{\'e}d{\'e}rale de Lausanne}
}

@article{henri2022reinit,
title = {Geometrical level set reinitialization using closest point method and kink detection for thin filaments, topology changes and two-phase flows},
journal = {Journal of Computational Physics},
volume = {448},
pages = {110704},
year = {2022},
issn = {0021-9991},
author = {Félix Henri and Mathieu Coquerelle and Pierre Lubin}
}

@article{zhao2005reinit,
 ISSN = {00255718, 10886842},
 author = {Hongkai Zhao},
 journal = {Mathematics of Computation},
 number = {250},
 pages = {603--627},
 publisher = {American Mathematical Society},
 title = {A Fast Sweeping Method for Eikonal Equations},
 urldate = {2025-11-09},
 volume = {74},
 year = {2005}
}

@article{sethian2001reinit,
title = {Evolution, Implementation, and Application of Level Set and Fast Marching Methods for Advancing Fronts},
journal = {Journal of Computational Physics},
volume = {169},
number = {2},
pages = {503-555},
year = {2001},
issn = {0021-9991},
author = {J.A. Sethian}
}

@article{chopp2001reinit,
author = {Chopp, David L.},
title = {Some Improvements of the Fast Marching Method},
journal = {SIAM Journal on Scientific Computing},
volume = {23},
number = {1},
pages = {230-244},
year = {2001}
}

@article{anumolu2013reinit,
title = {Gradient augmented reinitialization scheme for the level set method},
journal = {International Journal of Numerical Methods in Fluids},
author = {Lakshman Anumolu and Mario F. Trujillo},
volume = {73},
number = {12},
pages = {1011-1041},
year = {2013}
}

@article{shakoor2025reinit,
title = {Review of level-set reinitialization methods in computational mechanics and materials science},
journal = {Modelling and Simulation in Materials Science and Engineering},
author = {Modesar Shakoor},
volume = {33},
number = {5},
year = {2025}
}

@article{berger1984amr1,
  title={Adaptive mesh refinement for hyperbolic partial differential equations},
  author={Berger, Marsha J and Oliger, Joseph},
  journal={Journal of Computational Physics},
  volume={53},
  number={3},
  pages={484--512},
  year={1984},
  publisher={Elsevier}
}

@article{berger1989amr2,
  title={Local adaptive mesh refinement for shock hydrodynamics},
  author={Berger, Marsha J and Colella, Phillip},
  journal={Journal of Computational Physics},
  volume={82},
  number={1},
  pages={64--84},
  year={1989},
  publisher={Elsevier}
}

@article{zeng2023amr,
title = {A consistent adaptive level set framework for incompressible two-phase flows with high density ratios and high Reynolds numbers},
journal = {Journal of Computational Physics},
volume = {478},
pages = {111971},
year = {2023},
issn = {0021-9991},
author = {Yadong Zeng and Han Liu and Qiang Gao and Ann Almgren and Amneet Pal Singh Bhalla and Lian Shen}
}

@article{sussman1999amr,
  title={An adaptive level set approach for incompressible two-phase flows},
  author={Sussman, Mark and Almgren, Ann S and Bell, John B and Colella, Phillip and Howell, Louis H and Welcome, Michael L},
  journal={Journal of Computational Physics},
  volume={148},
  number={1},
  pages={81--124},
  year={1999},
  publisher={Elsevier}
}

@article{khokhlov1998amr,
title = {Fully Threaded Tree Algorithms for Adaptive Refinement Fluid Dynamics Simulations},
journal = {Journal of Computational Physics},
volume = {143},
number = {2},
pages = {519-543},
year = {1998},
issn = {0021-9991},
author = {A.M Khokhlov}
}

@article{strain1999amr,
title = {Tree Methods for Moving Interfaces},
journal = {Journal of Computational Physics},
volume = {151},
number = {2},
pages = {616-648},
year = {1999},
issn = {0021-9991},
author = {John Strain},
}

@article{mirzade2016amr,
title = {Parallel level-set methods on adaptive tree-based grids},
journal = {Journal of Computational Physics},
volume = {322},
pages = {345-364},
year = {2016},
issn = {0021-9991},
author = {Mohammad Mirzadeh and Arthur Guittet and Carsten Burstedde and Frederic Gibou}
}

@article{morgan2017tetra,
title = {{3D} level set methods for evolving fronts on tetrahedral meshes with adaptive mesh refinement},
journal = {Journal of Computational Physics},
volume = {336},
pages = {492-512},
year = {2017},
issn = {0021-9991},
author = {Nathaniel R. Morgan and Jacob I. Waltz}
}

@article{burstedde2011p4est,
author = {Burstedde, Carsten and Wilcox, Lucas C. and Ghattas, Omar},
title = {p4est: Scalable Algorithms for Parallel Adaptive Mesh Refinement on Forests of Octrees},
journal = {SIAM Journal on Scientific Computing},
volume = {33},
number = {3},
pages = {1103-1133},
year = {2011}

}

@article{min2010reinitializing,
  title={On reinitializing level set functions},
  author={Min, Chohong},
  journal={Journal of Computational Physics},
  volume={229},
  number={8},
  pages={2764--2772},
  year={2010},
  publisher={Elsevier}
}

@article{osher2004level,
  title={Level set methods and dynamic implicit surfaces},
  author={Osher, Stanley and Fedkiw, Ronald and Piechor, Krzysztof},
  journal={Appl. Mech. Rev.},
  volume={57},
  number={3},
  pages={B15--B15},
  year={2004}
}

@article{sussman2000coupled,
  title={A coupled level set and volume-of-fluid method for computing {3D} and axisymmetric incompressible two-phase flows},
  author={Sussman, Mark and Puckett, Elbridge Gerry},
  journal={Journal of Computational Physics},
  volume={162},
  number={2},
  pages={301--337},
  year={2000},
  publisher={Elsevier}
}

@article{Babuska1981,
  title={The p-version of the finite element method},
  author={Babuska, Ivo and Szabo, Barna A and Katz, I Norman},
  journal={SIAM journal on numerical analysis},
  volume={18},
  number={3},
  pages={515--545},
  year={1981},
  publisher={SIAM}
}

@article{ntoukas2022entropy,
  title={An entropy--stable p--adaptive nodal discontinuous Galerkin for the coupled {Navier--Stokes/Cahn--Hilliard} system},
  author={Ntoukas, Gerasimos and Manzanero, Juan and Rubio, Gonzalo and Valero, Eusebio and Ferrer, Esteban},
  journal={Journal of Computational Physics},
  volume={458},
  pages={111093},
  year={2022},
  publisher={Elsevier}
}

@article{mossier2023efficient,
  title={An efficient hp-adaptive strategy for a level-set ghost-fluid method},
  author={Mossier, Pascal and Appel, Daniel and Beck, Andrea D and Munz, Claus-Dieter},
  journal={Journal of Scientific Computing},
  volume={97},
  number={2},
  pages={50},
  year={2023},
  publisher={Springer}
}

@article{popinet2015quadtree,
  title={A quadtree-adaptive multigrid solver for the {Serre--Green--Naghdi} equations},
  author={Popinet, St{\'e}phane},
  journal={Journal of Computational Physics},
  volume={302},
  pages={336--358},
  year={2015},
  publisher={Elsevier}
}

@article{hergibo2024quadtree,
  title={A quadtree-based adaptive moment-of-fluid method for interface reconstruction with filaments},
  author={Hergibo, Philippe and Liang, Qiuhua and Phillips, Timothy N and Xie, Zhihua},
  journal={Journal of Computational Physics},
  volume={499},
  pages={112719},
  year={2024},
  publisher={Elsevier}
}

@article{yang2025octree,
  title={Octree-based adaptive mesh refinement and the shifted boundary method for efficient fluid dynamics simulations},
  author={Yang, Cheng-Hau and Scovazzi, Guglielmo and Krishnamurthy, Adarsh and Ganapathysubramanian, Baskar},
  journal={Advances in Computational Science and Engineering},
  volume={4},
  pages={57--84},
  year={2025},
  publisher={Advances in Computational Science and Engineering}
}

@article{ngo2020multi,
  title={A multi-level adaptive mesh refinement for an integrated finite element/level set formulation to simulate multiphase flows with surface tension},
  author={Ngo, Long Cu and Choi, Hyoung Gwon},
  journal={Computers \& Mathematics with Applications},
  volume={79},
  number={4},
  pages={908--933},
  year={2020},
  publisher={Elsevier}
}

@article{zhang2009parallel,
  title={A parallel algorithm for adaptive local refinement of tetrahedral meshes using bisection},
  author={Zhang, Lin-Bo and others},
  journal={Numer. Math.: Theory, Methods and Applications},
  volume={2},
  number={65-89},
  pages={13},
  year={2009}
}

@article{bargteil2006semi,
  title={A semi-Lagrangian contouring method for fluid simulation},
  author={Bargteil, Adam W and Goktekin, Tolga G and O'brien, James F and Strain, John A},
  journal={ACM Transactions on Graphics (TOG)},
  volume={25},
  number={1},
  pages={19--38},
  year={2006},
  publisher={ACM New York, NY, USA}
}

@article{kurioka2025improved,
  title={Improved particle level set method with higher-order kernel function correction: Enhancing accuracy and conservation},
  author={Kurioka, Shunsuke and Hu, Changhong},
  journal={Computers \& Fluids},
  volume={291},
  pages={106571},
  year={2025},
  publisher={Elsevier}
}

@article{gammanpila2025stabilized,
  title={Stabilized Nitsche-Type CIP/GP CutFEM for Two-Phase Flow Applications},
  author={Gammanpila, Himali and Aulisa, Eugenio and Chierici, Andrea},
  journal={Mathematics},
  volume={13},
  number={17},
  pages={2853},
  year={2025},
  publisher={MDPI}
}

@article{sampath2010parallel,
  title={A parallel geometric multigrid method for finite elements on octree meshes},
  author={Sampath, Rahul S and Biros, George},
  journal={SIAM Journal on Scientific Computing},
  volume={32},
  number={3},
  pages={1361--1392},
  year={2010},
  publisher={SIAM}
}

@article{gibou2018review,
  title={A review of level-set methods and some recent applications},
  author={Gibou, Frederic and Fedkiw, Ronald and Osher, Stanley},
  journal={Journal of Computational Physics},
  volume={353},
  pages={82--109},
  year={2018},
  publisher={Elsevier}
}

@article{ngo2017multi,
  title={A multi-level adaptive mesh refinement method for level set simulations of multiphase flow on unstructured meshes},
  author={Ngo, Long Cu and Choi, Hyoung Gwon},
  journal={International Journal for Numerical Methods in Engineering},
  volume={110},
  number={10},
  pages={947--971},
  year={2017},
  publisher={Wiley Online Library}
}

@article{Rivara1984,
  author  = {Rivara, M.-C.},
  title   = {Mesh refinement processes based on the generalized bisection of simplices},
  journal = {SIAM Journal on Numerical Analysis},
  volume  = {21},
  number  = {3},
  pages   = {604--613},
  year    = {1984}
}

@article{Kossaczky1994,
  author  = {Kossaczky, I.},
  title   = {A recursive approach to local mesh refinement in two and three dimensions},
  journal = {Journal of Computational and Applied Mathematics},
  volume  = {55},
  pages   = {275--288},
  year    = {1994}
}

@article{Bey2000,
  author  = {Bey, J.},
  title   = {Tetrahedral grid refinement},
  journal = {Computing},
  volume  = {55},
  pages   = {355--378},
  year    = {2000}
}

@misc{femus-github,
  author       = {Aulisa, Eugenio},
  title        = {{FEMuS}: A Finite Element Multiphysics Solver},
  howpublished = {\url{https://github.com/eaulisa/MyFEMuS}},
  year         = {2024},
  note         = {GitHub repository}
}

@article{aulisa2018monolithic,
  title={A monolithic ALE Newton--Krylov solver with multigrid-Richardson--Schwarz preconditioning for incompressible fluid-structure interaction},
  author={Aulisa, Eugenio and Bna, Simone and Bornia, Giorgio},
  journal={Computers \& Fluids},
  volume={174},
  pages={213--228},
  year={2018},
  publisher={Elsevier}
}

@article{calandrini2019fluid,
  title={Fluid-structure interaction simulations of venous valves: A monolithic ALE method for large structural displacements},
  author={Calandrini, Sara and Aulisa, Eugenio},
  journal={International Journal for Numerical Methods in Biomedical Engineering},
  volume={35},
  number={2},
  pages={e3156},
  year={2019},
  publisher={Wiley}
}

@article{calandrini2020field,
  title={A field-split preconditioning technique for fluid-structure interaction problems with applications in biomechanics},
  author={Calandrini, Sara and Aulisa, Eugenio and Ke, Guoyi},
  journal={International Journal for Numerical Methods in Biomedical Engineering},
  volume={36},
  number={3},
  pages={e3301},
  year={2020},
  publisher={Wiley}
}

@article{capodaglio2017particle, title={A particle tracking algorithm for parallel finite element applications}, author={Capodaglio, Giacomo and Aulisa, Eugenio}, journal={Computers \& Fluids}, volume={159}, pages={338--355}, year={2017}, publisher={Elsevier} }

@article{aulisa2019construction,
  title={Construction of h-refined continuous finite element spaces with arbitrary hanging node configurations and applications to multigrid algorithms},
  author={Aulisa, Eugenio and Capodaglio, Giacomo and Ke, Guoyi},
  journal={SIAM Journal on Scientific Computing},
  volume={41},
  number={1},
  pages={A480--A507},
  year={2019},
  publisher={SIAM}
}

@misc{blanco2014nanoflann,
  title   = {{nanoflann}: a {C}++ header-only fork of {FLANN}, a library for Nearest Neighbor ({NN}) search with {KD}-trees},
  author  = {Blanco, Jose Luis and Rai, Pranjal Kumar},
  year    = {2014},
  url     = {https://github.com/jlblancoc/nanoflann},
  urldate = {2026-02-04},
  note    = {GitHub repository}
}

@article{Wang2009,
  author  = {Zhaoyuan Wang and Jianming Yang and Bonguk Koo and Frederick Stern},
  title   = {A Coupled Level Set and Volume-of-Fluid Method for Sharp Interface Simulation of Plunging Breaking Waves},
  journal = {International Journal of Multiphase Flow},
  volume  = {35},
  number  = {3},
  pages   = {227--246},
  year    = {2009},
  doi     = {10.1016/j.ijmultiphaseflow.2008.11.004}
}

@article{Cervone2009,
  author  = {Antonio Cervone and Sandro Manservisi and Ruben Scardovelli and St{\'e}phane Zaleski},
  title   = {A Geometrical Predictor--Corrector Advection Scheme and Its Application to the Volume Fraction Function},
  journal = {Journal of Computational Physics},
  volume  = {228},
  number  = {2},
  pages   = {406--419},
  year    = {2009},
  doi     = {10.1016/j.jcp.2008.09.016}
}

@article{Ramanuj2019,
  author  = {Vimal A. Ramanuj and Ramanan Sankaran},
  title   = {High Order Anchoring and Reinitialization of Level Set Function for Simulating Interface Motion},
  journal = {Journal of Scientific Computing},
  volume  = {81},
  number  = {3},
  pages   = {1963--1986},
  year    = {2019},
  doi     = {10.1007/s10915-019-01053-6}
}

@article{Xia2023,
  author  = {Qing Xia and Junxiang Yang and Yibao Li},
  title   = {On the Conservative Phase-Field Method with the N-Component Incompressible Flows},
  journal = {Physics of Fluids},
  volume  = {35},
  number  = {1},
  year    = {2023},
  doi     = {10.1063/5.0135490}
}

\end{document}